\documentstyle[12pt,amsmath,amssymb]{article}

\begin{document}
\newtheorem{lem}{Lemma}[section]
\newtheorem{th}{Theorem}[section]
\newtheorem{prop}{Proposition}[section]
\newtheorem{rem}{Remark}[section]
\newtheorem{define}{Definition}[section]
\newtheorem{cor}{Corollary}[section]

\allowdisplaybreaks

\makeatletter\@addtoreset{equation}{section}\makeatother
\def\theequation{\arabic{section}.\arabic{equation}}
\newcommand{\ddGamma}{\overset{{.}{.}}{\Gamma}_X}
\newcommand{\dddGamma}{\overset{{.}{.}}{\Gamma}_{\R^d}}
\newcommand{\ddG}{\overset{{.}{.}}{\Gamma}_{\gamma,\gamma'}}
\newcommand{\PP}{{\bf P}}
\newcommand{\D}{{\cal D}}
\newcommand{\N}{{\Bbb N}}
\newcommand{\C}{{\Bbb C}}
\newcommand{\Z}{{\Bbb Z}}
\newcommand{\R}{{\Bbb R}}
\newcommand{\Rp}{{\R_+}}
\newcommand{\eps}{\varepsilon}
\newcommand{\Om}{\Omega}
\newcommand{\OM}{\Om_X^{\Rp}}
\newcommand{\om}{\omega}
\newcommand{\Pe}{{\bf P}}

\newcommand{\supp}{\operatorname{supp}}
\newcommand{\la}{\langle}
\newcommand{\ra}{\rangle}
\newcommand{\const}{\operatorname{const}}
\newcommand{\Crb}{C_{\rho,{\mathrm b}}(\dddGamma)}
\newcommand{\Xsym}{X^\infty_{\mathrm sym}}
\newcommand{\Gamman}{\Gamma_{X,0}}
\newcommand{\tob}{\Subset}
\newcommand{\myfrac}[2]{\mbox{$\frac{#1}{#2}$}}

\renewcommand{\P}{{\bf P}}
\renewcommand{\emptyset}{\varnothing}
\renewcommand{\tilde}{\widetilde}
\newcommand{\rom}[1]{{\rm #1}}
\newcommand{\FC}{{\cal F}C_{\mathrm b}^\infty\big({\cal D},\Gamma_X\big)}
\newcommand{\vph}{\varphi}
\newcommand{\fii}{\vph}
\newcommand{\di}{\partial}
\renewcommand{\div}{\operatorname{div}}
\newcommand{\pim}{{\pi_m}}
\newcommand{\HG}{H^\Gamma}
\newcommand{\HX}{H^X}
\newcommand{\Crbo}{C_{\rho,{\mathrm b}}(\overset{{.}{.}}{\Omega}_{X})}
\renewcommand{\author}[1]{\medskip{\large #1}\par\medskip}
\newcommand{\dd}{\overset{{.}{.}}}

$\text{}$\vspace{-28mm}

\begin{center}\LARGE \bf
The Heat Semigroup\\[2mm]\LARGE\bf on Configuration
Spaces\end{center}

\begin{center}\large Yuri Kondratiev${}^*$, Eugene Lytvynov${}^{**}$, Michael
R\"ockner${}^{***}$\end{center}

\begin{abstract}
In this paper, we study properties of the heat semigroup of
configuration space analysis. Using a natural ``Riemannian-like''
structure of the configuration space $\Gamma_X$ over a complete,
connected, oriented, and stochastically complete Riemannian
manifold $X$ of infinite volume, the heat semigroup
$(e^{-tH^\Gamma})_{t\in\R_+}$ was  introduced and studied in [{\it
J. Func.\ Anal.}\ {\bf 154} (1998), 444--500]. Here, $H^\Gamma$ is
the Dirichlet operator of the Dirichlet form ${\cal E}^\Gamma$
over the space $L^2(\Gamma_X,\pi_m)$, where $\pi_m$ is the Poisson
measure on $\Gamma_X$ with intensity $m$---the volume measure on
$X$. We construct a metric space $\Gamma_\infty$ that is
continuously embedded into $\Gamma_X$. Under some conditions on
the manifold $X$, we prove that $\Gamma_\infty$ is a set of full
$\pi_m$ measure and derive  an explicit  formula for the heat
semigroup:
$(e^{-tH^\Gamma}F)(\gamma)=\int_{\Gamma_\infty}F(\xi)\,\P_{t,\gamma}(d\xi)$,
where $\P_{t,\gamma}$ is a probability measure on $\Gamma_\infty$
for all $t>0$, $\gamma\in\Gamma_\infty$. The central results of
the paper are two types of  Feller properties for the heat
semigroup. The first one is a kind of  strong Feller property with
respect to the metric on the space $\Gamma_\infty$. The second
one, obtained in the case $X=\R^d$, is the Feller property with
respect to the intrinsic metric of the Dirichlet form ${\cal
E}^\Gamma$. Next, we give a direct construction of the independent
infinite particle process on the manifold $X$, which is a
realization of the Brownian motion on the configuration space. The
main point here is that we prove that this process can start in
every $\gamma\in\Gamma_\infty$, will never leave $\Gamma_\infty$,
and has continuous sample path in $\Gamma_\infty$, provided
$\operatorname{dim}X\ge2$. In this case,  we also prove that this
process is a strong  Markov process whose transition probabilities
are given by the $\P_{t,\gamma}(\cdot)$ above.  Furthermore, we
discuss the necessary changes to be done for constructing the
process in the case $\operatorname{dim}X=1$. Finally, as an easy
consequence we get a ``path-wise'' construction of the independent
particle process on $\Gamma_\infty$ from the underlying Brownian
motion.
\end{abstract}

\noindent \underline{\hspace{5cm}}\vspace{1mm}

\noindent Communicated by T. Kawai. Received November 8, 2001

\noindent 2000  Mathematics Subject Classification(s): Primary:
47D07, 60J60; Secondary: 60G57

\noindent${}^*$ Fakult\"at f\"ur Mathematik, Universit\"at
Bielefeld, Postfach 10 01 31, D-33501 Bielefeld, Germany;
 Institute of
Mathematics, Kiev, Ukraine; BiBoS, Univ.\ Bielefeld, Germany\\
{\it e-mail}: kondrat@mathematik.uni-bielefeld.de

\noindent ${}^{**}$ Institut f\"{u}r Angewandte Mathematik,
Universit\"{a}t Bonn, Wegelerstr.~6, D-53115 Bonn, Germany;
 BiBoS, Univ.\ Bielefeld, Germany\\ {\it e-mail}:
 lytvynov@wiener.iam.uni-bonn.de

 \noindent ${}^{***}$ Fakult\"at f\"ur Mathematik, Universit\"at
Bielefeld, Postfach 10 01 31, D-33501 Bielefeld, Germany; BiBoS,
Univ.\ Bielefeld, Germany\\ {\it e-mail}:
roeckner@mathematik.uni-bielefeld.de

\newpage \noindent{\it Running head}: Heat semigroup on configuration
spaces

\tableofcontents

\newpage

\section{Introduction} In  \cite{AKR1, AKR2, AKR3, AKR4}, stochastic
analysis and differential geometry on  configuration spaces were
considerably developed by using the so-called ``lifting
procedure,'' see also  \cite{Roeckner, MaR, Ro, RS, ADL1, ADL2}
for further results and reviews.

Let us recall that the configuration space $\Gamma_X$ over a
complete, connected, oriented, and stochastically complete
Riemannian manifold $X$ of infinite volume is defined as the set
of all infinite subsets of $X$ which are locally finite. Each
configuration $\gamma\in\Gamma_X$ can be identified with the Radon
measure $\sum_{x\in\gamma}\eps_x$. The tangent space  to
$\Gamma_X$ at a point $\gamma\in\Gamma_X$, denoted by
$T_\gamma(\Gamma_X)$, is defined as  the direct sum of the tangent
spaces  to $X$ at $x$, where $x$ runs over the points of the
configuration $\gamma$; that is,
$T_\gamma(\Gamma_X):=\bigoplus_{x\in\gamma}T_x(X)$. The gradient
$\nabla^\Gamma F(\gamma)$ of a differentiable function $F\colon
\Gamma_X\to\R$ at a point $\gamma\in\Gamma_X$ is defined as an
element of the tangent space $T_\gamma(\Gamma_X)$ through a
natural lifting of  the gradient on $X$.  Analogously, one
introduces also the notion of divergence of a vector field over
$\Gamma_X$.

Let $\pi_m$ denote the Poisson measure on $\Gamma_X$ with
intensity $m$---the volume measure on $X$.  By using the
integration by parts formula for the Poisson measure, it was shown
in \cite{AKR3} that $\pi_m$ is a volume measure on $\Gamma_X$, in
the sense that the gradient and the divergence become dual
operators on $L^2(\pi_m):=L^2(\Gamma_X,\pi_m)$.

Thus, having  identified differentiation and a volume measure on
the configuration space, the next step in \cite{AKR3} was to
consider the Dirichlet form over $L^2(\pi_m)$, which is defined by
$${\cal E
}^\Gamma(F_1,F_2):=\frac12\,\int_{\Gamma_X}\langle\nabla^\Gamma
F_1(\gamma),\nabla^\Gamma
F_2(\gamma)\rangle_{T_\gamma(\Gamma_X)}\,\pi_m(d\gamma)$$ on an
appropriate set of smooth cylinder functions on $\Gamma_X$. Using
again the integration by parts formula, one obtains the associated
Dirichlet operator, i.e., the operator  $H^\Gamma$ in $L^2(\pi_m)$
satisfying ${\cal E}^\Gamma(F_1,F_2)=(H^\Gamma
F_1,F_2)_{L^2(\pi_m)}$. This yields, in particular, that the
bilinear form ${\cal E}^\Gamma$ is closable. Moreover, the
operator $H^\Gamma$ was shown to be  essentially selfadjoint. We
will preserve the notation $H^\Gamma$ for its closure.

The present paper is devoted to the study of properties of the
heat semigroup $(e^{-tH^\Gamma})_{t\in\R_+}$.

By using the general theory of Dirichlet forms, it has been
already proved in \cite{AKR3, MaR}  (see also  \cite{RS}) that
there exists a diffusion process (i.e.\ a strong Markov process
with continuous sample paths) on the configuration space that is
canonically associated with the heat semigroup
$(e^{-tH^\Gamma})_{t\in\R_+}$, i.e., for each $F\in L^2(\pi_m)$,
$$(e^{-tH^\Gamma}F)(\gamma)=\int_{\pmb\Omega}F({\bf X}_t)\,d{\bf
P}_\gamma$$ for $\pi_m$-a.a.\ (or even quasi-every)
$\gamma\in\Gamma_X$. This process is then the Brownian motion on
$\Gamma_X$. Moreover, this process is, in fact, the well-known
independent infinite particle process (cf.\ \cite{AKR3}). The
latter  is obtained by taking countably many independent Brownian
motions on $X$, see \cite{Doob}.

The first part of this paper is devoted to deriving an explicit
formula  for the heat semigroup. We introduce functionals $B_n$,
$n\in\N$, on $\Gamma_X$ by $$B_n(\gamma):=\sum_{x\in\gamma}
\exp\big[-\myfrac1n\, \operatorname{dist}(x_0,x)\,\big], $$ where
$x_0$ is a fixed point of the manifold $X$. We define a subset
$\Gamma_\infty$ of $\Gamma_X$ consisting of those configurations
$\gamma$ for which $B_n(\gamma)<\infty$ for all $n\in\N$, and
equip $\Gamma_\infty$ with a metric in such a way that the
convergence in $\Gamma_\infty$ means vague convergence together
with convergence of all the  functionals $B_n$ (see also
\cite{KS}). Under some conditions on the geometry of the manifold
$X$, we prove that $\Gamma_\infty$ is of full $\pi_m$ measure and
that, for each $\gamma\in\Gamma_\infty$, $t>0$, there exists a
probability measure ${\bf P}_{t,\gamma}$ on $ \Gamma_\infty$ such
that  for each $F\in L^2(\Gamma_\infty,\pi_m)$
\begin{equation}
\label{iowa}(e^{-tH^\Gamma}F)(\gamma)=\int_{\Gamma_\infty}F(\xi)\,\P_{t,\gamma}(d\xi),\qquad
\text{$\pi_m$-a.a.\ $\gamma\in\Gamma_\infty$.} \end{equation}

To this end, we apply the method  for  constructing   probability
measures on $\Gamma_X$ described in \cite{VGG}, and define ${\bf
P}_{t,\gamma}$ via a product measure $\bigotimes_{k=1}^\infty
p_{t,x_k}$ on $X^\N$. Here, $p_{t,x}(dy):=p(t,x,y)\, m(dy)$,
$p(t,x,y)$ is the heat kernel of the manifold $X$,  and
$\gamma=\{x_k\}_{k=1}^\infty$ (the resulting measure
$\P_{t,\gamma}$ will,  however, be  independent of the chosen
ordering of the points of $\gamma$).

The second part of the paper is devoted to our main results which
concern two types of  Feller properties  of the heat semigroup.

 We introduce a class $\bf D$
of measurable functions on $\Gamma_\infty$, which particularly
contains
 all bounded local functions, and  show that $\bf D$
is invariant under the action of the semigroup $e^{-tH^\Gamma}$.
Moreover, we prove that, for each $F\in{\bf D }$, the map
$$\Gamma_\infty\ni\gamma\mapsto
(\P_tF)(\gamma):=\int_{\Gamma_\infty}F(\xi)\,
\P_{t,\gamma}(d\xi)\in\R$$ is a continuous function on the space
$\Gamma_\infty$ ($\P_tF$ is even continuous with respect to some
  weaker metric). Thus, we obtain a kind of  strong Feller property
of the heat semigroup. Here, we use  results on harmonic analysis
over configuration spaces from \cite{Kuna1, Kuna11, Kuna3} (see
also \cite{Len1,Len2,Len3,Macchi, BKKL}),

 Next, we consider a
metric space $\dd\Gamma_\infty$ which is an appropriate extension
of $\Gamma_\infty$ to multiple configurations in $X$, i.e., to
$\Z_+$-valued Radon measures on $X$.  Restricting ourselves to the
case $X=\R^d$, we prove that the operators $({\bf P}_t)_{t>0}$
defined by $$\dd\Gamma_\infty\ni\gamma\mapsto({\bf
P}_tF)(\gamma):=\int_{\dd\Gamma_{\infty}}F(\xi)\,{\bf
P}_{t,\gamma}(d\xi)\in \R$$ preserve the class of all bounded
functions on $\dd\Gamma_{\infty}$ which are continuous with
respect to the intrinsic metric of the Dirichlet form ${\cal
E}^\Gamma$ (see \cite{Ro}). Thus, for this metric, we have  the
usual Feller property of the heat semigroup.

In the third part of the paper, which is more probabilistic, we
present a direct construction of the independent infinite particle
process on the manifold $X$  with the state space
$\dd\Gamma_\infty$, which will be therefore a realization of the
Brownian motion on the configuration space mentioned above. We
show that, if the dimension of $X$ is $\ge2$, the constructed
process is the unique continuous strong Markov process on
$\Gamma_\infty$ whose transition probabilities are given by
$\P_t(\gamma,\cdot):=\P_{t,\gamma}(\cdot)$. In particular, it
starts at any configuration in $\Gamma_\infty$ and never leaves
$\Gamma_\infty$.  If $\operatorname{dim}X=1$, one cannot exclude
collisions of the particles, but it is still possible to realize
the Brownian motion on the configuration space as a continuous
Markov process on $\dd\Gamma_\infty$. Finally, we describe a
path-wise construction of the infinite particle process starting
from any point in $\Gamma_\infty$ (respectively
$\dd\Gamma_\infty$). More precisely, we show that the obvious
heuristic construction can be performed rigorously.

We should mention that independent infinite particle processes
have been studied by many authors, see e.g.\ \cite{shiga}, but
neither in this paper, nor in any other reference we are aware of,
it was proved that the process takes values in the configuration
space for all values of $t>0$.

\section{Intrinsic Dirichlet operator on the Poisson space}
\label{section2}

In this section, we will briefly recall the definition and some
properties of the intrinsic Dirichlet operator on the Poisson
space. We refer the reader to \cite{AKR3, AKR1, AKR2} for details
and  proofs.

Let $X$ be a complete, connected, oriented $C^\infty$ Riemannian
manifold.  Let $m$ denote the volume measure on $X$, and we
suppose that $m(X)=\infty$. Let $\nabla^X$ and
$H^X:=-\frac12\Delta^X$ be the gradient and Laplace--Beltrami
operator on $X$, respectively. We denote by $\D:=C_0^\infty(X)$
the space of all $C^\infty$ functions on $X$ with compact support.
It is well-known that $(H^X,\D)$ is essentially selfadjoint on
$L^2(m):=L^2(X,{\cal B}(X),m)$, where ${\cal B}(X)$ is the Borel
$\sigma$-algebra on $X$. In what follows, we will always suppose
that $H^X$ is conservative (cf.\ e.g.\ \cite{T}).

Let $p(t,x,y)$, $t\in(0,\infty)$, $x,y\in X$, denote the heat
kernel of the operator $H^X$:
\begin{equation}\label{2}
(e^{-tH^X}\fii)(x)=\int_X\fii(y)p(t,x,y)\,m(dy),\qquad
\text{$m$-a.e.}\ x\in X,\end{equation} where $\fii$ is a bounded
measurable function on $X$.  We recall that $p(t,x,y)$ is a
strictly positive $C^\infty$ function on $(0,\infty)\times X\times
X$ (cf.\ e.g.\ \cite{Davies}).

The conservativity condition yields, in particular, that
\begin{equation}\label{6745635}\int_X p(t,x,y)\,m(dy)=1,\qquad t\in(0,\infty),\, x\in X,
\end{equation}
i.e., for each $t>0$ and $x\in X$ the heat kernel determines a
probability measure
\begin{equation}\label{plplp}p_{t,x}(dy):=p_t(x,dy)=p(t,x,y) \,m(dy)\end{equation} on $X$.
Thus, the manifold $X$ is stochastically complete.

Next, we  consider the configuration space $\Gamma_X$ over
$X$---the set of all  infinite  subsets in $X$ which are locally
finite: $$\Gamma_X:=\{\gamma\subset X\mid |\gamma|=\infty \text{
and }|\gamma_\Lambda|<\infty\text{ for each compact
}\Lambda\subset X\}.$$ Here, $|\cdot|$ denotes the cardinality of
a set and $\gamma_\Lambda:= \gamma\cap\Lambda$. One can identify
any $\gamma\in\Gamma_X$ with the positive Radon measure
$$\sum_{x\in\gamma}\eps_x\in{\cal M}(X), $$ where  ${\cal M}(X)$
 stands for the set of all positive
 Radon  measures  on
${\cal B}(X)$. The space $\Gamma_X$ can be endowed with the
relative topology as a subset of the space ${\cal M}(X)$ with the
vague topology, i.e., the weakest topology on $\Gamma_X$ with
respect to which  all maps
$$\Gamma_X\ni\gamma\mapsto\la\fii,\gamma\ra:=\int_X\fii(x)\,\gamma(dx)
=\sum_{x\in\gamma}\fii(x),\qquad\fii\in{\cal D},$$ are continuous.
We shall denote  the Borel $\sigma$-algebra on $\Gamma_X$ by
${\cal B}(\Gamma_X)$.

Let $\pim$ denote the Poisson measure on $(\Gamma_X,{\cal
B}(\Gamma_X))$ with intensity $m$. This measure can be
characterized by its Laplace transform
\begin{equation}\label{ewrwerewrwe5}\ell_\pim(\fii):=\int_{\Gamma_X}
e^{\la\fii,\gamma\ra}\,\pim(d\gamma)
=\exp\bigg(\int_X(e^{\fii(x)}-1)\,m(dx)\bigg),\qquad
\fii\in\D.\end{equation} We refer to e.g.\ \cite{VGG,Sh,AKR3} for
a detailed discussion of the construction of the Poisson measure
on the configuration space. Now, we recall  how to define the
intrinsic Dirichlet operator $H^\Gamma$ in the space
$L^2(\pim):=L^2(\Gamma_X,{\cal B}(\Gamma_X),\pi_m)$.

Let $T_x(X)$ denote the tangent space to $X$ at a point $x\in X$.
The tangent space to $\Gamma_X$ at a point $\gamma\in\Gamma_X$ is
defined as the Hilbert space
$$T_\gamma(\Gamma_X):=\bigoplus_{x\in\gamma}T_x(X).$$
Thus, each $V(\gamma )\in T_\gamma (\Gamma _X)$ has the form $%
V(\gamma )=(V(\gamma,x ))_{x\in \gamma }$, where $V(\gamma,x )\in
T_x(X)$, and $$ \| V(\gamma )\| _{T_\gamma(\Gamma_X)}
^2=\sum_{x\in \gamma }\|V(\gamma,x )\|_{T_x(X)}^2.$$

Let $\gamma \in \Gamma _X$ and $x\in \gamma $. We denote by ${\cal
O}_{\gamma ,x}$  an arbitrary open neighborhood of $x$ in $X$ such
that ${\cal O}_{\gamma,x}\cap(\gamma \setminus
\{x\})=\varnothing$.  Now, for a function $F\colon\Gamma _X\to
{\Bbb R}$, $\gamma \in \Gamma _X$, and $x\in\gamma $, we define a
function $F_{x}(\gamma ,\cdot )\colon{\cal O}_{\gamma ,x}\to {\Bbb
R}$ by $${\cal O}_{\gamma,x}\ni y\mapsto
F_x(\gamma,y):=F((\gamma\setminus\{x\})\cup\{y\})\in{\Bbb R}.$$

We say that a function $F\colon \Gamma _X\to {\Bbb R}$ is
differentiable at $\gamma \in \Gamma _X$ if for each $x\in \gamma
$ the function $F_x(\gamma ,\cdot )$ is differentiable at $x$ and
\[
\nabla ^\Gamma F(\gamma ):=(\nabla ^X F_x(\gamma,x ))_{x\in \gamma
}\in T_\gamma( \Gamma _X),
\]
where
\begin{equation}\notag
\nabla ^X F_x(\gamma,x ):=\nabla ^X_yF_x(\gamma ,y){\big |}_{y=x}
\end{equation}
(cf.\ \cite{ADL1,ADL2}). Evidently, this definition is independent
of the choice of the set ${\cal O}_{\gamma,x}$. We will call
$\nabla^\Gamma F(\gamma)$ the gradient of $F$ at $\gamma$.

 We introduce
the set $\FC$ consisting of all smooth cylinder functions on
$\Gamma_X$, i.e., all functions  of the form
\begin{equation}\label{3}
F(\gamma)=g_F(\la\fii_1,\gamma\ra,\dots,\la\fii_N,\gamma\ra),
\qquad\gamma\in\Gamma_X,
\end{equation}
where $N\in\N$, $\fii_1,\dots,\fii_N\in\D$, and $g_F\in
C^\infty_{\mathrm b}(\R^N)$. Any function $F$ of the form \eqref3
is differentiable at each point $\gamma\in\Gamma_X$, and its
gradient is given by
\begin{equation}\label{4}
(\nabla^\Gamma F)(\gamma)=\sum_{j=1}^N \di_j\, g_F(\la
\fii_1,\gamma\ra,\dots,\la\fii_N,\gamma\ra)\nabla^X\fii_j,
\end{equation}
where $\di_j\, g_F$ means derivative with respect to the $j$-th
coordinate.

Then, the corresponding  pre-Dirichlet form is
\begin{equation}\label{three}
{\cal E}^\Gamma(F,G):=\frac12\int_{\Gamma_X}\la\nabla^\Gamma
F(\gamma), \nabla^\Gamma G(\gamma)
\ra_{T_\gamma(\Gamma_X)}\,\pim(d\gamma),\qquad
F,G\in\FC.\end{equation} By using the integration by parts formula
on the Poisson space, one shows that the associated Dirichlet
operator $H^\Gamma$, i.e., the operator satisfying $${\cal
E}^\Gamma(F,G)=(H^\Gamma F,G) _{L^2(\pi_m)},\qquad F,G\in\FC,$$ is
of the form
\begin{gather*}
(H^\Gamma F)(\gamma)=-\sum_{i,j=1}^N \di_i\,\di_j\, g_F(\la
\fii_1,\gamma\ra,\dots,\la\fii_N,\gamma\ra)\int_X\frac12\,\la\nabla^X\fii_i(x),
\nabla^X\fii_j(x)\ra_{T_x(X)}\,\gamma(dx)\\ \text{}+\sum_{j=1}^N
\di_j \, g_F (\la
\fii_1,\gamma\ra,\dots,\la\fii_N,\gamma\ra)\int_X(H^X\fii_j)(x)\,\gamma(dx),
\end{gather*}
where $F$ is given by \eqref3. Therefore, the bilinear form
$({\cal E}^\Gamma,\FC)$ is closable on $L^2(\pi_m)$, and with its
closure we can associate  a positive definite selfadjoint
operator, the Friedrichs extension of $H^\Gamma$, which
 will be also denoted by $H^\Gamma$. (In fact, $\FC$ is a domain of essential
 selfadjointness of $H^\Gamma$, see \cite[Theorem~5.3]{AKR3}.)

Consider the corresponding heat semigroup $(e^{-tH^\Gamma})_{t\in
\R_+}$ in $L^2(\pim)$, where as usual $\R_+:=[0,\infty)$. We set
$$E(\D_1,\Gamma):=\operatorname{l.h.}\big
\{\,\exp\big[\,\la\log(1+\fii),\cdot\ra\,\big]\mid
\fii\in\D_1\,\big\}.$$ Here, l.h.\ means linear hull and
\begin{equation}\label{kljhgf}\D_1:=\big\{\, \fii\in D(H^X)\cap L^1(m)\mid H^X\fii\in
L^1(m)\ \text{and $-\delta\le\fii\le0$ for some }\delta\in(0,1)
\,\big\}. \end{equation}

\begin{prop}\label{prop1}
We have
\begin{equation*}
e^{-tH^\Gamma}\exp\big[\,\la\log(1+\fii),\cdot\ra\,\big]=
\exp\big[\,\la\log(1+e^{-t\HX}\fii),\cdot\ra\,\big]\qquad
\text{\rm $\pi_m$-a.e.\ for all }\fii\in{\cal D}_1.
\end{equation*}
\end{prop}

\noindent{\it Proof}. See \rom{\cite[Proposition~4.1]{AKR3}}.\quad
$\blacksquare$\vspace{2mm}

As a direct consequence of this proposition,  in particular, one
obtains that\linebreak $(\HG,E(\D_1,\Gamma))$ is essentially
selfadjoint on $L^2(\pim)$.

Finally,  the diffusion process that is properly associated with
the Dirichlet form $({\cal E}^\Gamma,\operatorname{Dom}({\cal
E}^\Gamma))$ is the usual independent infinite particle process,
or in other terms,  Brownian motion on $\Gamma_X$ (cf.\
\cite[Subsection~6.2]{AKR3}).

\section[Correlation measures in
configuration space analysis]{Correlation measures in
configuration space\\ analysis}\label{correlation}

In this section, we shall recall some facts on $K$-transforms and
correlation measures. We shall follow  \cite{Kuna1, Kuna3} (see
also \cite{Len1, Len2, Len3,Macchi, Kuna11, Kuna2, BKKL}; in
\cite{Kuna1, Kuna3, BKKL}
 the reader can also find  many further references and
historical comments).

Denote by $\Gamman$ the space of all finite configurations over
$X$: $$\Gamman:=\bigsqcup_{n=0}^\infty \Gamma_X^{(n)},\qquad
\Gamma_X^{(0)} =\{\varnothing\},\qquad
\Gamma_X^{(n)}=\{\eta\subset X\mid |\eta|=n\},\qquad n\in\N.$$

Let $$ \widetilde X ^{n}=\big\{\,(x_1,\dots,x_n)\in X^n \mid
x_i\ne x_j\text{ when } i\ne j\,\big\}$$ and let $S_n$ denote the
group of permutations of $\{1,\dots,n\}$, which acts on $\tilde
X^n$ by
$$\sigma(x_1,\dots,x_n)=(x_{\sigma(1)},\dots,x_{\sigma(n)}),\qquad
\sigma\in S_n .$$ Through the natural bijection
\begin{equation}\label{101}
\widetilde X^{n}/S _n\mapsto \Gamma_X^{(n)}\end{equation} one
defines a topology on $\Gamma_X^{(n)}$. The space $\Gamman$ is
then  equipped with the topology of disjoint unions. Let ${\cal B
}(\Gamman)$ denote the corresponding Borel $\sigma$-algebra. A set
$K\subset \Gamman$ is compact iff there exists $N\in\N$ with
$K\cap \Gamma_X^{(n)}=\varnothing$ for all $n> N$ and $K\cap
\Gamma_X^{(n)}$ is compact for all $n\le N$.
 The set of all Borel
sets in $\Gamman$ with compact closure  is denoted by ${\cal B}
_{\mathrm c}(\Gamman)$.

A ${\cal B}(\Gamman)$-measurable function $G\colon\Gamman\to\R$ is
said to have bounded support if there exist a relatively compact
open  set $\Lambda\subset X$ and $N\in\N$ such that $\{G\ne 0
\}\subset \bigsqcup_{n=0}^N \Gamma_\Lambda^ {(n)}$. The space of
 bounded functions on $\Gamman$ with bounded support is
denoted  by $B_{\mathrm bs}(\Gamman)$.

In what follows, for any $\gamma\in\Gamma_X$, we shall  use the
notation $\sum_{\eta\tob\gamma}$ for the summation over all
$\eta\subset\gamma$ such that $|\eta|< \infty$. For a function
$G\colon \Gamman\to \R$, the $K$-transform of $G$ is then  defined
by
\begin{equation}\label{102}
(KG)(\gamma):=\sum_{\eta\tob\gamma} G(\eta)
\end{equation}
for each $\gamma\in\Gamma_X$ such that at least one of the series
$\sum_{\eta\tob\gamma} G^+(\eta)$ or $\sum_{\eta\tob\gamma}
G^-(\eta)$ converges in $\R_+$, where
$G^+(\gamma):=\max\{0,G(\gamma)\}$,
$G^-(\gamma)=-\min\{0,G(\gamma)\}$. For each $G\in B_{\mathrm
bs}(\Gamman)$ and  each $\gamma\in\Gamma_X$, the series
$\sum_{\eta\tob\gamma} G(\eta)$ is always finite, and moreover,
$(KG)(\cdot)$ is a ${\cal B}(\Gamma_X)$-measurable function on
$\Gamma_X$ (cf.~\cite[Proposition~3.5]{Kuna1}).

Let $\mu$ be a probability measure on $(\Gamma_X,{\cal
B}(\Gamma_X))$. The correlation measure corresponding to $\mu$ is
defined by $$\rho_\mu(A):=\int_{\Gamma_X}(K{\bf
1}_A)(\gamma)\,\mu(d\gamma),\qquad A\in{\cal B}(\Gamma_{X,0}).$$
 $\rho_\mu$ is obviously a  measure on
$(\Gamman,{\cal B} (\Gamman))$.

\begin{prop}\label{prop3.1} Let $\mu$ be a probability measure on
$(\Gamma_X,{\cal B}(\Gamma_X))$\rom. Then\rom, the measure
$\rho_\mu$ is locally finite\rom, i\rom.e\rom{.,}
\begin{equation}\label{104}
\rho_\mu(A)<\infty\quad \text{\rom{for all }}A\in{\cal B}_{\mathrm
c}(\Gamman),\end{equation} if and only if
\begin{equation}\label{lambada}
\int_{\Gamma_X}|\gamma_\Lambda|^n\, \mu(d\gamma)<\infty
\quad\text{\rom{for all $n\in\N$ and $\Lambda\in{\cal B}_{\mathrm
c} (X)$}}.
\end{equation}
\end{prop}

\noindent {\it Proof}. See \cite[Proposition~4.2]{Kuna1}. \quad
$\blacksquare$\vspace{2mm}

We say that a measure $\mu$ satisfying \eqref{lambada}  has finite
local moments  and denote the set of all such measures on
$(\Gamma_X,{\cal B}(\Gamma_X))$  by ${\cal M}_{\mathrm fm
}(\Gamma_X)$. The set of all locally finite measures on $\Gamman$
will be denoted by ${\cal M}_{\mathrm lf}(\Gamman)$.

\begin{prop}\label{prop3.2}
Let $\mu\in{\cal M}_{\mathrm fm}(\Gamma_X)$ and let
$G\colon\Gamman\to\R$ be a measurable function which  is
integrable with respect to the measure $\rho_\mu$\rom. Then\rom,
$KG$ is
 well-defined and finite $\mu$-a\rom.e\rom.\rom, and integrable with respect
 to the measure $\mu$\rom. If for some
$G'\colon\Gamman\to\R$\rom, $G=G'$ $\rho_\mu$-a\rom.e\rom.\rom,
then $KG=KG'$ $\mu$-a\rom.e\rom.\rom, and hence the $K$-transform
defines a linear mapping $$K\colon L^1(\Gamman,{\cal
B}(\Gamman),\rho_\mu)\to L^1(\Gamma_X,{\cal B}(\Gamma_X),\mu).$$
Furthermore, we have
\begin{gather*}
\|KG\|_{L^1(\mu)}\le\|K|G|\,\|_{L^1(\mu)}=\|G\|_{L^1(\rho_\mu)}
\end{gather*}
and
$$\int_{\Gamman}G(\eta)\,\rho_\mu(d\eta)=\int_{\Gamma_X}(KG)(\gamma)\,
\mu(d\gamma).$$
\end{prop}

\noindent{\it Proof}. See
\cite[Theorem~4.11]{Kuna1}.\quad$\blacksquare$\vspace{2mm}

For two functions $G_1,G_2\colon\Gamman\to\R$, the
$\star$-convolution
 of $G_1$ and $G_2$ is defined as the mapping $G_1\star G_2\colon\Gamman\to\R$ given by
\begin{equation}\label{starproduct}(G_1\star G_2)(\eta):=\sum_{(\eta_1,\eta_2,\eta_3)\in{\cal
P}_3(\eta)}
G_1(\eta_1\cup\eta_2)G_2(\eta_2\cup\eta_3),\end{equation}
 where
${\cal P}_3(\eta)$ denotes the set of all ordered partitions
$(\eta_1,\eta_2,\eta_3)$ of $\eta$ into 3 parts. Clearly, if
$G_1$, $G_2$ are  ${\cal B}(\Gamma_{X,0})$-measurable, then so is
 $G_1\star G_2$. The main property of the
$\star$-convolution is given by the following formula (see
\cite[Proposition~3.11]{Kuna1}):
\begin{equation}\label{waerfed}
(K(G_1\star G_2))(\gamma)=(KG_1)(\gamma)\cdot (KG_2)(\gamma),
\end{equation} provided $(KG_1)(\gamma)$ and $(KG_2)(\gamma)$
exist.

Let $\sigma$ be a non-atomic Radon measure. The
Lebes\-gue--Poisson measure $\lambda_\sigma$ on $(\Gamman,{\cal
B}(\Gamman))$ with intensity $\sigma$ is defined by
$$\lambda_\sigma:=\eps_\varnothing+\sum_{n=1}^\infty\frac1{n!}\,\sigma^{\otimes
n },$$ where the measure $\sigma^{\otimes n}$ is defined on
$\Gamma_X^{(n)}$ via the bijection \eqref{101}.

Finally, let us introduce the notion of correlation functions.
Suppose that a measure $\rho\in{\cal M}_{\mathrm lf}(\Gamman)$ is
absolutely continuous with respect to the Lebesgue--Poisson
measure $\lambda_m$ with intensity $m$, and define the functions
$k^{(n)}\colon\Gamma_X^{(n)}\to\R$ as the restrictions of the
Radon--Nikodym derivative $k:=\dfrac{d\rho}{d\lambda_m}$ to
$\Gamma_X^{(n)}$. In the case where $\rho=\rho_\mu$ is a
correlation measure, the functions $(k_\mu^{(n)})_{n=1}^\infty$
are called correlation functions of the measure $\mu$.

\section{Heat kernel measures $\P_{t,\gamma}$}

In this section, we shall  construct a family of probability
measures $\P_{t,\gamma}$ on the configuration space so that
$\P_{t,\bullet}(\cdot)$ is the kernel of the integral operator
$e^{-tH^\Gamma}$.

First, we recall the construction of  probability measures on the
configuration space $\Gamma_X$ proposed by A.~M.~Vershik et~al.\
\cite{VGG}, see also \cite{HS}.

Let us consider the infinite product
$X^\N=\overset{\infty}{\underset {k=1}{\times}}X_k$, $X_k=X$,
furnished with the product topology, and let ${\cal B}(X^\N)$
denote the Borel $\sigma$-algebra on $X^\N$. We define $\tilde
X^\N$ as the set of all elements $(x_1,x_2,\dots)\in X^\N$ such
that 1) $x_i\ne x_j$ when $i\ne j$, and 2) the sequence
$\{x_k\}_{k=1}^\infty$
 has no accumulation points in $X$.
Evidently,  \begin{multline}\label{esrdgzer}\tilde
X^\N=\left[\bigcap_{i\ne j}\big\{\,(x_1,x_2,\dots)\in X^\N: x_i\ne
x_j\,\big\}\right]\cap\\ \cap\left[ \bigcap_{n=1}^\infty
\bigcup_{k=1}^\infty \big\{\,(x_1,x_2,\dots)\in X^\N: \forall l\ge
k\ d(x_0,x_l)\ge n\,\big\}\right],\end{multline} where $x_0$ is a
fixed point of $X$ and $d(\cdot,\cdot)$ denotes the distance on
$X$. Hence,  $\tilde X^\N\in{\cal B}(X^\N)$.

Let $\nu_k$, $k\in\N$, be  nonatomic probability measures on
$(X,{\cal B}(X))$ and consider the  product measure
$\nu:=\bigotimes_{k=1}^\infty\nu_k$  on $(X^\N,{\cal B}(X^\N))$.
 \eqref{esrdgzer} and the Borel--Cantelli lemma imply the
 following:

\begin{lem}
\label{lem1} \rom{\cite{VGG}} We have $\nu(\tilde X^\N)=0$ or
$1$\rom, and $\nu(\tilde X^\N)=1$ if and only if
\begin{equation}\label{qwer}\sum\limits_{k=1}^\infty\nu_k(\Lambda)<\infty\qquad
\text{\rom{for each compact} }
 \Lambda\subset X.\end{equation}
\end{lem}

Let $S_\infty$ denote the group of all permutations of the
sequence of natural numbers, which  acts on $X^\N$:
$$\sigma(y_1,y_2,\dots)=(y_{\sigma(1)},y_{\sigma(2)},\dots),\qquad
\sigma\in S_\infty.$$ The space $\tilde X^\N$ is invariant under
the action of $S_\infty$. Through the natural bijection $\tilde
X^\N/S_\infty\mapsto\Gamma_X$, we shall identify these two spaces.
Let $I\colon \tilde X^\N\to\Gamma_X$ be given by
\begin{equation}\label{oizgdr}\tilde X^\N\ni{\bf x}=(x_1,x_2,\dots)\mapsto
I{\bf x}=\{x_1,x_2,\dots\}\in\Gamma_X.\end{equation} Thus, $I$
maps an element ${\bf x}\in \tilde X^\N$ into the corresponding
equivalence class $[ {\bf x}] \in \tilde X^\N/S_\infty$.

The mapping $I\colon \tilde X^\N\to\Gamma_X$ defined by
\eqref{oizgdr} is ${\cal B}(\tilde X^\N)$-${\cal
B}(\Gamma_X)$-measurable (here ${\cal B}(\tilde X^\N)$ denotes the
trace $\sigma$-algebra of ${\cal B}(X^\N)$ on $\tilde X^\N $).
Indeed, the $\sigma$-algebra ${\cal B}(\Gamma_X)$ is generated by
the sets of the form
$$A_{\Lambda,n}=\big\{\,\gamma\in\Gamma_X\colon \langle {\bf 1
}_\Lambda,\gamma\rangle=n\,\big\},$$ where $n\in{\Bbb
Z}_+:=\N\cup\{0\}$, $\Lambda$ is a compactum in $X$, and ${\bf
1}_\Lambda$ is the indicator of $\Lambda$ (see e.g.\
\cite{Kal86,MKM}). Then,
\begin{equation}\label{gfsdsdghfv}I^{-1}(A_{\Lambda,n})=\bigcup_{\sigma\in S^{\mathrm fin}
_\infty}\big\{\, (x_1,x_2,\dots)\in\tilde X^\N\colon
x_{\sigma(i)}\in\Lambda,\ i=1,\dots,n,\
x_{\sigma(i)}\in\Lambda^{\mathrm c},\ i\ge n+1
\,\big\},\end{equation} where $S^{\mathrm fin}_\infty$ denotes the
group of all finite permutations of the sequence of natural
numbers, and $\Lambda^{\mathrm c}:=X\setminus \Lambda$. Since
$S^{\mathrm fin}_\infty$ has a countable number of elements, we
conclude from \eqref{gfsdsdghfv} that
$I^{-1}(A_{\Lambda,n})\in{\cal B}(\tilde X^\N)$, which implies the
measurability of $I$.

Hence, if the measures $\nu_k$, $k\in\N$, satisfy  condition
\eqref{qwer}, we can consider the image of the probability measure
$\nu$ on $\tilde X^\N$ under the mapping $I$, which is a
probability measure on $\Gamma_X$. Evidently, this image-measure
is independent of the order of the $\nu_k$'s, that is, it
coincides with the measure on $\Gamma_X$ constructed through the
product-measure $\bigotimes_{k=1}^\infty \nu_{\sigma(k)}$ for each
$\sigma\in S_\infty$.

Let now $t>0$ and let $\gamma$ be a fixed point of $\Gamma_X$ such
that
\begin{equation}\label{6}
\sum_{x\in\gamma} p_{t,x}(\Lambda)<\infty\quad\text{for each
compact } \Lambda\subset X,\end{equation} where $p_{t,x}$ is as in
\eqref{plplp}. Define
\begin{equation}\label{kltzwe}
\text{${\bf P}_{t,\gamma}:={\Bbb P}_{t,{\bf x}}\circ I^{-1}$,
where }{\Bbb P}_{t,{\bf x}}:= \bigotimes_{k=1}^\infty
p_{t,x_k}\end{equation}  and where ${\bf
x}=(x_k)_{k=1}^\infty\in\tilde X^\N$ is an arbitrary element of
the set $I^{-1}\{\gamma\}$ (the resulting measure $\P_{t,\gamma}$
being independent of the choice of $\bf x$).

Below, we shall need the correlation measure $\rho_{t,\gamma}$ of
$\P_{t,\gamma}$.

\begin{prop} \label{awse} Let $t>0$ and let $\gamma\in\Gamma_X$ satisfy
\eqref{6}\rom. Then, the correlation measure $\rho_{t,\gamma}$ of
$\P_{t,\gamma}$ is given by
\begin{align}
\rho_{t,\gamma}\restriction \Gamma_X^{(0)}:&=\rho_{t,\gamma}^{(0)}
:=\varepsilon_{\varnothing},\notag\\ \rho_{t,\gamma}\restriction
\Gamma_X^{(n)}:&=\rho_{t,\gamma}^{(n)}
:=\sum_{\theta\tob\gamma:\,|\theta|=n}\left(\underset{x\in\theta}{\otimes}
\,p_{t,x}\right)\circ T_n^{-1}, \qquad
n\in\N,\label{chyk}\end{align} where 
$T_n\colon \tilde X^n\to \Gamma^{(n)}_X$ is the composition of the
natural quotient map $\tilde X^n\to \tilde X^n/S_n$ and the
bijection \eqref{101} \rom(the measure
$(\otimes_{x\in\theta}p_{t,x})\circ T_n^{-1}$ is independent of a
chosen order of the product of the measures $p_{t,x}$\rom)\rom.
Moreover\rom, we have $$\rho_{t,\gamma}(\Gamma_\Lambda^{(n)})\le
\frac1{n!}\,\left(\sum_{x\in\gamma}p_{t,x}(\Lambda)\right)^n<\infty\qquad\text{\rom{for
each compact} }\Lambda\subset X,$$ in particular\rom,
$\rho_{t,\gamma} \in{\cal M}_{\mathrm lf}(\Gamman)$ and
$\P_{t,\gamma}\in{\cal M}_{\mathrm fm}(\Gamma_X).$\end{prop}

\noindent {\it Proof}. Let the measure $\rho_{t,\gamma}$ on
$(\Gamma_{X,0},{\cal B}(\Gamman))$  be defined by \eqref{chyk}.
For a measurable function $G\colon\Gamman\to\R$, we have
$G=(G^{(n)})_{n=0}^\infty$, where
$G^{(n)}:=G\restriction\Gamma_X^{(n)}$. Then, by using the
definition of $\P_{t,\gamma}$ and  the monotone convergence
theorem, we have, for any ${\cal B}(\Gamma_{X,0})$-measurable
function $G\colon\Gamma_{X,0}\to\R_+$, that
\begin{gather}
\int_{\Gamma_X}(KG)(\xi)\,\P_{t,\gamma}(d\xi)=\int_{\tilde
X^\N}(KG)(I{\bf y})\,{\Bbb P}_{t,{\bf x}}(d{\bf y})\notag\\
=G^{(0)}(\varnothing)+\sum_{n=1}^\infty
\sum_{\{i_1,\dots,i_n\}\subset\{1,2,\dots\}}\int_{\tilde X^\N}
G^{(n)}\circ T_n(y_{i_1},\dots,y_{i_n})\,\bigotimes_{k=1}^\infty
p_{t,x_k}(dy_k)\notag\\=G^{(0)}(\varnothing)+ \sum_{n=1}^\infty
\sum_{\{i_1,\dots,i_n\}\subset\{1,2,\dots\}} \int_{\tilde
X^n}G^{(n)}\circ T_n(y_1,\dots,y_n)\,
p_{t,x_{i_1}}\otimes\cdots\otimes
p_{t,x_{i_n}}(dy_1,\dots,dy_n)\notag\\=G^{(0)}(\varnothing)+
\sum_{n=1}^\infty \int_{\tilde X^n} G^{(n)}\circ T_n
(y_1,\dots,y_n)\,\sum_{\{i_1,\dots,i_n\}\subset\{1,2,\dots\}}
\,p_{t,x_{i_1}}\otimes \cdots\otimes
p_{t,x_{i_n}}(dy_1,\dots,dy_n)\notag\\=\int_{\Gamman}G(\eta)\,\rho_{t,\gamma}(d\eta),
\label{awdrawdr}\end{gather} where ${\bf x}=(x_k)_{k=1}^ \infty\in
I^{-1}\{\gamma\}$. The final inequality in the assertion
immediately follows from   \eqref{6} and \eqref{chyk}. Hence, the
measure $\rho_{t,\gamma}$ is
 from ${\cal M}_{\mathrm lf}(\Gamman)$, and therefore, by
Proposition~\ref{prop3.1}, $\P_{t,\gamma}\in{\cal M}_{\mathrm
fm}(\Gamma_X)$.\quad$\blacksquare$

\begin{rem}\rom{ One could  also start with a measure
$\rho_{t,\gamma}$ on $\Gamman$ that is given  {\it a priori} by
formula \eqref{chyk} for each $t>0$ and each  $\gamma\in\Gamma_X$
satisfying \eqref6, and then, using \cite[Theorem~6.5]{Kuna1},
identify $\P_{t,\gamma}$ as the unique probability measure on
$\Gamma_X$ whose correlation measure is $\rho_{t,\gamma}$.}

\end{rem}

Our next aim is to show that  condition \eqref6 is satisfied for
$\pi_m$-a.a.\ $\gamma\in\Gamma_X$, at least under some additional
conditions on the manifold $X$.

Let us assume that the manifold $X$ satisfies the following two
conditions:

\begin{description}

\item[\rom{(C1)}] For each $t>0$, there exist constants $C_t>0$
and $\eps_t>0$ such that
\begin{equation*}
p(t,x,y)\le C_t\exp\big[-d(x,y)^{1+\eps_t}\,\big],\qquad
 t>0,\ x,y\in X.
\end{equation*}

\item[\rom{(C2)}] For some fixed $x_0\in X$,
$$m\big(B(x_0,r)\big)\le c_{x_0}\, r^N,\qquad r>0,$$ where
$c_{x_0}>0$, $N\in\N$, and $B(x_0,r)$ denotes the geodesic ball
with center at $x$ and  radius $r$.

\end{description}

Concerning these conditions, in particular, the upper estimate
 of the heat kernel, we refer the reader e.g. to
\cite{Davies,Gr1,Gr2} and the references therein. For example, in
the case of a manifold $X$ of nonnegative Ricci curvature, one has
\begin{gather}
p(t,x,y)\le\frac{C}{m\big(B(y,\sqrt
t)\big)}\exp\bigg(-\frac{d(x,y)^2}{(4+\eps)t} \bigg),\qquad
\eps>0,\label{seghdr}\\ m\big(B(x,r)\big)\le
\operatorname{const}_d r^d \label{scuople}
\end{gather}
($d$ being the dimension of $X$). Thus,  conditions (C1) and (C2)
are satisfied if the manifold $X$ possesses  the following
additional property:
\begin{equation}\label{waeergdzt}\forall r>0:\qquad
\inf_{x\in X} m\big(B(x,r)\big)>0,\end{equation} which is true,
for example, in case of a manifold having bounded geometry (see
\cite{Davies}).

Now, we shall follow the idea of  \cite{KS} to consider
 subsets of the configuration space on which some special
 functionals take finite values.
  So, for each $n\in\N$, we
 introduce the functional
\begin{equation}\label{judr}
B_n(\gamma):=\big\langle\exp\big[-\mbox{$\frac{1}{n}$}\,d(x_0,\cdot)\,\big],\gamma
\big\rangle=\sum_{x\in\gamma}\exp\big[
-\mbox{$\frac{1}n$}\,d(x_0,x)\, \big],\qquad
\gamma\in\Gamma_X,\end{equation} and define $\Gamma_n\in{\cal
B}(\Gamma_X)$ by
\begin{equation}\label{wemkfwe}
\Gamma_n:=\big\{\,\gamma\in\Gamma_X:\,B_n(\gamma)<\infty\,\big\}.
\end{equation} Here, $x_0$ is as in (C2).
Evidently, we have, for each $n\in\N$,
$\Gamma_{n+1}\subset\Gamma_n$, and let
$$\Gamma_\infty:=\bigcap_{n=1}^\infty \Gamma_n.$$

Let $d_{\mathrm V}$ be any metric on ${\cal M}(X)$ determining the
vague topology. For example, we can take as $d_{\mathrm V}$ the
metric $d_{\mathrm K }$ that was introduced in \cite{Rachev}:
$$d_{\mathrm K}(\nu_1,\nu_2):=\sum_{i=1}^\infty
2^{-i}\,d_{{\mathrm K},i}(\nu_1,\nu_2)/[1+d_{{\mathrm
K},i}(\nu_1,\nu_2)],\qquad \nu_1,\nu_2\in{\cal M}(X),$$ where
\begin{multline*} d_{{\mathrm K},i}(\nu_1,\nu_2):=\sup\bigg\{\,
\Big|\int_X f\, d(\nu_1-\nu_2)\Big|\colon\, f:X\to\R,\\
\sup_{x,y\in X }\frac{d(f(x),f(y))}{d(x,y)}\le1,\  f(x)=0\
\text{if }d(x_0,x)\ge i\,\bigg\}.\end{multline*} The metric
$d_{\mathrm K}$ is a generalization of the Kantorovich metric, and
on any set of measures from ${\cal M}(X)$ which have  uniformly
bounded support, $d_{\mathrm K}$ is just equivalent to the
Kantorovich metric.

Then, we can metrize the set $\Gamma_\infty$ as follows: for
$\gamma_1,\gamma_2\in\Gamma_\infty$
\begin{equation}\label{secondmetric}d_{\infty}(\gamma_1,\gamma_2):=d_ {\mathrm
V}(\gamma_1,\gamma_2)+\sum_{n=1}^\infty 2^{-n}
|B_n(\gamma_1)-B_n(\gamma_2)|\big/\big[1+|B_n(\gamma_1)-B_n(\gamma_2)|\,\big].\end{equation}

Let ${\cal B}(\Gamma_\infty)$ denote the trace $\sigma$-algebra of
${\cal B}(\Gamma_X)$ on $\Gamma_\infty$. It can be shown that this
$\sigma$-algebra coincides with the Borel $\sigma$-algebra on
$\Gamma_\infty$ that corresponds to the topology generated by the
$d_\infty$ metric.


\begin{prop}\label{quadrofonia} Let \rom{(C1)} and \rom{(C2)} be
satisfied\rom. Then\rom, $\Gamma_\infty$ is a set of full $\pi_m$
measure\rom. Furthermore\rom, for each $\gamma\in\Gamma_\infty$
condition \eqref{6} is satisfied and  $\Gamma_\infty$ is a set of
full $\P_{t,\gamma}$ measure  for each $t>0$.\end{prop}

\noindent{\it Proof}. We have by (C2) that
\begin{align}\int_{\Gamma_X}
B_n(\gamma)\,\pi_m(d\gamma)&=\int_X\exp\big[-\myfrac1n\,d(x_0,x)\,\big]\,m(dx)\notag\\
&=\sum_{k=1}^\infty \int\limits_{B(x_0,k)\setminus
B(x_0,k-1)}\exp\big[-\myfrac1n\,d(x_0,x)\,\big]\,m(dx)\notag\\
&\le\sum_{k=1}^\infty\exp\big[-\myfrac1n\,(k-1)\,\big]m\big(B(x_0,k)\big)\notag\\
&\le\sum_{k=1}^\infty\exp\big[-\myfrac1n\,(k-1)\,\big]c_{x_0}k^N<\infty.\label{lpseaq}
\end{align}
Therefore, $B_n$ is $\pi_m$-a.e.\ finite, i.e.,
$\pi_m(\Gamma_n)=1$ for all $n\in\N$, which yields that
$\pi_m(\Gamma_\infty)=1$.

Next, from (C1) we get, for each $r>0$, $t>0$, and
$\gamma\in\Gamma_\infty$,
\begin{align*}\sum_{x\in\gamma}
p_{t,x}(B(x_0,r))&=\sum_{x\in\gamma}\int_{B(x_0,r)}p(t,x,y)\,m(dy)\\&\le
\tilde
C_t\sum_{x\in\gamma}\int_{B(x_0,r)}\exp\big[-d(x,y)\,\big]\,m(dy)
\\&\le\tilde  C_t\sum_{x\in\gamma}\exp\big[-d(x_0,
x)\,\big]\,\int_{B(x_0,r)}\exp\big[d(x_0,y)\,\big]\,m(dy)<\infty,\end{align*}
so that \eqref6 is satisfied.

Finally, for each $\gamma\in\Gamma_\infty$, $t>0$, and $n\in\N$,
we get from  (C1),  (C2) (cf.\ also \eqref{lpseaq}), and the
monotone convergence theorem that
\begin{align}\int_{\Gamma_X}B_n(\xi)\,\P_{t,\gamma}(d\xi)&=\sum_{x\in\gamma}
\int_X\exp\big[-\myfrac1n\,d (x_0,y)\,\big]p(t,x,y)\,m(dy)\notag\\
&\le\sum_{x\in\gamma}\int_X \exp\big[-\myfrac1n\,d
(x_0,y)\,\big]C_t\exp\big[-d(x,y)^{1+\eps_t}\,\big]\,m(dy)\notag\\
&\le C_{t,n}\sum_{x\in\gamma}\int_X\exp\big[-\myfrac1n\,d
(x_0,y)\,\big]\exp\big[-\myfrac1{2n}\,d
(x,y)\,\big]\,m(dy)\notag\\ &\le
C_{t,n}\sum_{x\in\gamma}\exp\big[-\myfrac1{2n}\,d(x_0,x)\,\big]\int_X\exp
\big[-\myfrac1{2n}\,d
(x_0,y)\,\big]\,m(dy)<\infty,\label{iikkkf}\end{align} which
yields that $\P_{t,\gamma}(\Gamma_\infty)=1$.\quad$\blacksquare$

\begin{rem}\rom{ Let $\pi_{zm}$ denote the Poisson measure on $(\Gamma_X,{\cal B}
(\Gamma_X))$ with intensity $zm$, where $z>0$. Since the
correlation measure of $\pi_{zm}$ is the Lebesgue--Poisson measure
$\lambda_{zm}$, it follows from the proof of
Proposition~\ref{quadrofonia} that $\pi_{zm}(\Gamma_\infty)=1$ for
all $z>0$. Furthermore, let $\mu_{\nu,m}:=\int_{0}^\infty
\pi_{zm}\,\nu(dz)$ be a mixed Poisson measure such that $\nu$ is a
probability measure on $(0,\infty)$. Then, $\Gamma_\infty$ is a
set of full $\mu_{\nu,m}$ measure. }\end{rem}

\section{Explicit formula for the heat semigroup}
Due to Proposition~ \ref{quadrofonia}, we can consider $\pi_m$ as
a probability measure on $(\Gamma_\infty,{\cal
B}(\Gamma_\infty))$. In this section, we shall derive  an explicit
formula for  the heat semigroup $(e^{-tH^\Gamma})_{t\in\R_+}$.


\begin{th}
\label{th1} Let the conditions  \rom{(C1)} and \rom{(C2)} be
satisfied\rom. Then\rom, for each $F\in
L^2(\Gamma_\infty,\pim)$\rom, we have
\begin{equation}\label{for1}
(e^{-tH^\Gamma}F)(\gamma) =\int_{\Gamma_\infty}F(\xi)\,{\bf
P}_{t,\gamma}(d\xi)\end{equation} for $\pim$-a.a.\
$\gamma\in\Gamma_\infty$\rom.
\end{th}

\noindent {\it Proof}. We start with the following

\begin{lem}\label{lem3}
Let $\tilde{{\cal D}}_1$ denote the subset of ${\cal D}_1$
\rom(see \eqref{kljhgf}\rom) given by $$\tilde{{\cal
D}}_1:=\big\{\, \fii\in{\cal D}\mid \exists \delta\in(0,1):
-\delta\le\fii\le0 \,\big\}.$$ Then\rom, for any
$\fii\in\tilde{{\cal D}}_1$\rom, $\gamma\in\Gamma_\infty$\rom, and
$t>0$\rom, we have
$$\int_{\Gamma_\infty}\exp\big[\,\la\log(1+\fii),\xi\ra\,\big]\,{\bf
P}_ {t,\gamma}(d\xi) = \exp\left[\,\left\la\log\left(1+ \int
\varphi\, dp_{t,\bullet} \right),\gamma\right\ra\,\right].$$
\end{lem}

\noindent {\it Proof}. First, we observe that by
Proposition~\ref{quadrofonia}
$$\sum_{x\in\gamma}\int_X|\varphi(y)|\,p_{t,x}(dy)<\infty\qquad
\text{for each $\gamma\in\Gamma_\infty$ and $\varphi\in\D$}.$$ But
for each $\gamma\in\Gamma_\infty$, ${\bf x}=(x_k)_{k=1}^\infty\in
I^{-1}\{\gamma\} $, and $\varphi\in\tilde\D_1$,
\begin{align*}\int_{\Gamma_\infty}\exp\big[\la\log(1+\varphi),\xi\ra\big]\,
\P_{t,\gamma}(d\xi)&= \prod_{k=1}^\infty \int_X (1+\fii(y))\,
p_{t,x_k}(dy)\\&= \exp\left[\,\left\la\log\left(1+ \int \varphi\,
dp_{t,\bullet} \right),\gamma\right\ra\,\right]. \quad\blacksquare
\end{align*}

\begin{lem}\label{waedawseaw} For any measurable function $F\colon
\Gamma_\infty\to\R_+$\rom, we have
\begin{equation}\label{erse}\int_{\Gamma_\infty}\int_{\Gamma_\infty}
F(\xi)\, \P_{t,\gamma}(d\xi)\,
\pi_m(d\gamma)=\int_{\Gamma_\infty}F(\gamma)\,\pi_m(d\gamma).\end{equation}
\end{lem}

\noindent {\it Proof}. It is easy to check that $\big\{\,
\exp\big[\la \log(1+\fii),\cdot\ra\big]\mid \fii \in\tilde {\cal
D}_1 \,\big\}$ is stable under multiplication and that it contains
a countable subset separating the points of $\Gamma_\infty$, so it
generates ${\cal B}(\Gamma_\infty)$. Therefore, we only have to
check \eqref{erse} for $F:=\exp\big[\la
\log(1+\fii),\cdot\ra\big]$, $\fii \in\tilde {\cal D}_1$.  But for
such functions \eqref{erse} immediately follows from
Lemma~\ref{lem3}. Indeed, \eqref{ewrwerewrwe5}  extends to all
functions $\fii\colon X\to\R_+$ which are increasing limits of
functions $\fii_n\in{\cal D}$, $n\in\N$, such as
$\log(1+\int\fii\,dp_{t,\bullet})$. Furthermore,
$\int\!\!\int\fii\, dp_{t,\bullet}\,dm=\int\fii\,dm$, since $H^X$
is assumed to be conservative.\quad$\blacksquare$ \vspace{2mm}

Now, we can easily finish the proof  of the theorem. It follows
from Lemma~\ref{waedawseaw} that, if $A\in{\cal B}(\Gamma_\infty)$
is of zero $\pi_m$ measure, then $\P_{t,\gamma}(A)=0$ for
$\pi_m$-a.e.\ $\gamma\in\Gamma_\infty$. Moreover, using the
Cauchy-Schwarz inequality and Lemma~\ref{waedawseaw}, we get
\begin{align*}\int_{\Gamma_\infty}\bigg(\int_{\Gamma_\infty}F(\xi)\,\P_{t,\gamma}(d\xi)\bigg)^2
\pi_m(d\gamma)&\le\int_{\Gamma_\infty}\int_{\Gamma_\infty}|F(\xi)|^2\,\P_{t,\gamma}(d\xi)
\,\pi_m(d\gamma)\\&=\int_{\Gamma_\infty}|F(\gamma)|^2\,\pi_m(d\gamma).\end{align*}

Thus, for each $t>0$, we can define a linear continuous operator
$$\P_t:L^2(\Gamma_\infty,\pi_m)\to L^2(\Gamma_\infty,\pi_m)$$ by
setting
$$(\P_tF)(\gamma):=\int_{\Gamma_\infty}F(\xi)\,\P_{t,\gamma}(d\xi).$$
By Proposition~\ref{prop1} and  Lemma~\ref{lem3}, the action of
the operator $\P_t$ coincides with the action of the operator
$e^{-tH^\Gamma}$ on the set
$\big\{\,\exp\big[\la\log(1+\fii),\cdot\ra\big]\mid\fii\in\tilde
{\cal D}_1\,\big\}$, which is total in $L^2(\Gamma_\infty,\pi_m)$
(i.e., its linear hull  is a dense set in
$L^2(\Gamma_\infty,\pi_m)$). Hence, we get the equality
$e^{-tH^\Gamma}=\P_t$, which proves the
theorem.\quad$\blacksquare$\vspace{2mm}

In what follows, for a measurable function $F$ on $\Gamma_\infty$,
we set
\begin{equation}\label{fuka}
({\bf P}_tF)(\gamma):=\int_{\Gamma_\infty}F(\xi)\,{\bf
P}_{t,\gamma}(d\xi), \qquad t>0,\,
\gamma\in\Gamma_\infty,\end{equation} provided the integral on the
right hand side exists. Hence, by virtue of Theorem~\ref{th1},
${\bf P}_tF$ is a $\pi_m$-version of $e^{-tH^\Gamma} F$ for each
$F\in L^2(\Gamma_\infty,\pi_m)$.

\begin{rem}\rom{One can easily prove an explicit formula for
the heat semigroup $(e^{-tH^\Gamma})_{t\in\R_+}$ in the weak
sense. More specifically, we define for each $t>0$ a function
$R_t\colon\Gamman\times\Gamman\to\R$ setting:
$R_t(\eta,\theta):=0$ if $|\eta|\ne|\theta|$\rom,
$R_t(\{\varnothing\},\{\varnothing\})=1$\rom, and for
$\eta=\{x_1,\dots,x_n\}$\rom, $\theta=\{y_1,\dots,y_n\}$,
$n\in\N$\rom, $$R_t(\eta,\theta):=\sum_{\sigma\in
S_n}\prod_{k=1}^n p_t(x_k,y_{ \sigma(k)}),$$ where
$p_t(x,y):=p(t,x,y)$\rom. Suppose that  conditions \rom{(C1)} and
\rom{(C2)} are satisfied\rom. Then\rom, for arbitrary measurable
functions $G_1,G_2\colon \Gamman \to \R$ such that
$$\int_{\Gamman}\int_{\Gamman}|G_1(\theta)|\cdot
\big[R_t(\eta,\theta)\star_\eta
|G_2(\eta)|\big]\,\lambda_m(d\eta)\,\lambda_m(d\theta)<\infty
\qquad\text{\rom{for all} }t>0$$ \rom($\star_\eta$ denoting the
$\star$-convolution with respect to the $\eta$ variable\rom{),} we
have $$ \int_{ \Gamma_\infty}|{\bf
P}_tF_1(\gamma)|\,|F_2(\gamma)|\,\pi_m(d\gamma)<\infty,$$ where
$F_1(\gamma)=(KG_1)(\gamma)$\rom,
$F_2(\gamma)=(KG_2)(\gamma)$\rom, and $$\int_{\Gamma_\infty}({\bf
P}_t F_1)(\gamma)F_2(\gamma)\,\pi_m(d\gamma)
=\int_{\Gamman}\int_{\Gamman} G_1(\theta)\cdot\big[
R_t(\eta,\theta) \star_\eta
G_2(\eta)\big]\,\lambda_m(d\eta)\,\lambda_m(d\theta).$$ }\end{rem}

Now, we define a family of probability kernels
$(\PP_t)_{t\in\R_+}$ on the space $(\Gamma_\infty,{\cal
B}(\Gamma_\infty))$  setting
\begin{equation}\label{label1}\PP_t(\gamma,A):=\P_{t,\gamma}(A),\qquad
\gamma\in\Gamma_\infty,\,A\in{\cal
B}(\Gamma_\infty),\,t\in\R_+,\end{equation} where
\begin{equation}\label{label2}\P_{0,\gamma}:=\varepsilon_\gamma.\end{equation}
Since $\gamma\mapsto {\bf P}_tF(\gamma)$ is measurable for $F$ in
the linear span of
$\big\{\,\exp\big[\la\log(1+\fii),\cdot\ra\big]\mid\fii\in\tilde{\cal
D} _1\,\big\}$ by Lemma~\ref{lem3}, a monotone class argument
shows that, indeed, $\gamma\mapsto{\bf P}_t(\gamma,A)$ is ${\cal
B}(\Gamma_\infty)$-measurable for all $A\in{\cal
B}(\Gamma_\infty)$.

We finish this section with the following proposition.

\begin{prop}\label{wadr} Let \rom{(C1)} and \rom{(C2)}
be satisfied\rom. Then\rom, $(\PP_t)_{t\in\R_+}$ is a Markov
semigroup of kernels on $(\Gamma_\infty,{\cal
B}(\Gamma_\infty))$\rom.\end{prop}

\noindent{\it Proof}.  The Markov property of  the kernels
$\PP_t$, i.e., $\P_{t}(\gamma,\Gamma_\infty)=1$,
$\gamma\in\Gamma_\infty$, follows from
Proposition~\ref{quadrofonia}.

Let us show the semigroup property: $\PP_t\PP_s=\PP_{t+s}$,
$t,s\in\R_+$. To this end, we fix $t,s>0$,
$\gamma\in\Gamma_\infty$, and $A\in{\cal B}(\Gamma_\infty)$. Then,
by  the construction of the measure $\P_{t,\gamma}$ and
 the
 semigroup property of the heat  kernel on
 $X$, we get
\begin{align*}(\P_t\P_s)(\gamma,A)&=
\int_{\Gamma_\infty}\P_{s,\xi}(A)\,\P_{t,\gamma}(d\xi)
=\int_{\tilde X^\N}\P_{s,I{\bf y}}(A)\,{\Bbb P}_{t,{\bf x}}(d{\bf
y})\\&=\int_{\tilde X^\N}{\Bbb P}_{s,{\bf y}}(I^{-1}A)\,{\Bbb
P}_{t,{\bf x }}(d{\bf y})=\int_{ X^\N}{\Bbb P}_{s,{\bf
y}}(I^{-1}A)\,{\Bbb P}_{t,{\bf x }}(d{\bf y})\\&={\Bbb
P}_{t+s,{\bf
x}}(I^{-1}A)=\P_{t+s,\gamma}(A)=\P_{t+s}(\gamma,A),\end{align*}
where ${\bf x}\in
I^{-1}\{\gamma\}$.\quad$\blacksquare$\vspace{2mm}

\section{A strong Feller property of the heat semigroup}

Let us introduce a new metric $d_1$ on the set $\Gamma_\infty$ as
follows: $$d_1(\gamma_1,\gamma_2):=d_{\mathrm
V}(\gamma_1,\gamma_2)+|B_1(\gamma_1)-B_1(\gamma_2)|,\qquad
\gamma_1,\gamma_2\in\Gamma_\infty.$$
 Evidently,  convergence with
respect to the $d_\infty$ metric implies  convergence with respect
to the $d_1$ metric.

In this section,  we shall show that the ``concrete version''
$({\bf P}_t)_{t\in\R_+}$ of the heat semigroup
$(e^{-tH^\Gamma})_{t\in\R_+}$ constructed in the previous section
possesses a kind of strong Feller property with respect to the
 metric $d_1$, and therefore also with respect to $d_\infty$.

\begin{th}\label{Meller_Skorohod}
Let \rom{(C1)} and \rom{(C2)} hold\rom. Let $G\colon \Gamman\to\R$
be  a  measurable function
 satisfying the following condition\rom:
\begin{equation}\label{705}\forall c>0:\quad
\int_{\Gamman}|G(\eta)|\,\lambda_{m_c}(d\eta)<\infty,
\end{equation} where\rom, for each $c>0$\rom, $\lambda_{m_c}$ is
the Lebesgue--Poisson measure on $\Gamman$ with intensity
\begin{equation}\label{rdtrdt}m_c(dx):=c\,e^{d(x_0,x)}\,m(dx).\end{equation}
 Then\rom, for each $t>0$,
the function
\begin{equation}\label{yqa} \Gamma_{\infty}\ni\gamma\mapsto ({\bf
P}_t(KG))(\gamma)=\int_{\Gamma_{\infty}}(KG)(\xi)\, {\bf
P}_{t,\gamma}(d\xi)\in\R\end{equation} is continuous with respect
to the metric $d_{1}$\rom.
\end{th}

\noindent {\it Proof}. Let $\gamma\in\Gamma_\infty$. By
Propositions~\ref{awse}, \ref{quadrofonia} and the definition of
correlation functions, we see that for each $t>0$  the measure
$\rho_{t,\gamma}$ is absolutely continuous  with respect to the
Lebesgue--Poisson measure $\lambda_m$, and the correlation
functions $k^{(n)}_{t,\gamma}$ of $\P_{t,\gamma}$ are given by
\begin{equation}\label{kj}
k^{(n)}_{t,\gamma}(\theta)=
\sum_{(i_1,\dots,i_n)\in\tilde{\N}{}^n}\,\prod _{k=1}^n
p_t(x_{i_k},y_k)\quad \text{\rom{for $\gamma=\{x_i\}_{i=1}^\infty$
and $\theta=\{y_1,\dots,y_n\}$}},
\end{equation}
where $$\tilde{\N}{}^n:=\big\{\,(i_i,\dots,i_n)\in\N^n: i_k\ne
i_l\text{ if }k\ne l\,\big\}.$$ Denote
\begin{equation}\label{lkjhj}k_t(\gamma,\theta):=k_{t,\gamma}(\theta):
=\frac{d\rho_{t,\gamma}}{d\lambda_m}(\theta),\end{equation} so
that $k_{t,\gamma}(\theta)=k^{(n)}_{t,\gamma}(\theta)$ for
$|\theta|=n$.

By using (C1), we get
\begin{align}
|k_t(\gamma,\{y_1,\dots,y_n\}) |&\le
\prod_{k=1}^n\bigg(\sum_{x\in\gamma} p_{t}(x,y_k)\bigg)\notag\\
&\le\prod_{k=1}^n\bigg(\sum_{x\in\gamma}
C_t'\exp\big[-d(x,y_k)\,\big]\bigg)\notag\\
&\le\bigg(C_t'\sum_{x\in\gamma}
\exp\big[-d(x_0,x)\,\big]\bigg)^n\,\exp\big[\,d(x_0,y_1)+\dots+d(x_0,y_n)\,\big].
\label{esrdrttud}\end{align} Hence,  \eqref{705} and
\eqref{rdtrdt} imply that $G\in L^1(\Gamma_{X,0},{\cal
B}(\Gamma_{X,0}),\rho_{t,\gamma})$. Therefore, if
$\gamma_j\to\gamma$ in $\Gamma_\infty$ with respect to $d_1$, by
Proposition~\ref{prop3.2} and \eqref{lkjhj} we have to prove that
\begin{equation}\label{708}
\int_{\Gamman} G(\eta)k_{t}(\gamma^j,\eta) \,\lambda_m(d\eta)\to
\int_{\Gamman} G(\eta)k_{t}(\gamma,\eta)
\,\lambda_m(d\eta)\qquad\text{as }j\to\infty.\end{equation}

First, we  show that
\begin{equation}\label{709}
k_t(\gamma^j,\eta)\to k_t(\gamma,\eta) \quad\text{as $j\to\infty$
for each fixed $\eta\in\Gamman$}.
\end{equation}

Since $\gamma^j\to\gamma$ in the $d_{1}$ metric, we have,
particularly, that $\gamma^j\to\gamma$ in the $d_{\mathrm V}$
metric. We claim that  that there exists  a numeration of the
points of the configurations
 $\gamma^j$, $j\in\N$,  and of $\gamma$  such  that
\begin{equation}\label{710}
\gamma^j=\{x_k^j\}_{k=1}^\infty,\quad
\gamma=\{x_k\}_{k=1}^\infty,\qquad \forall k\in\N:\quad
d(x_k^j,x_k)\to0\text{ as $j\to\infty$.}
\end{equation}

Indeed, let us fix any numeration of points of $\gamma$ such that
$$\gamma=\{x_k\}_{k=1}^\infty, \qquad d(x_0,x_{k+1})\ge
d(x_0,x_k),\quad k\in\N.$$ Next, we fix positive numbers $r_n$,
$n\in\N$, so that \begin{gather*}r_{n+1}>r_n,\quad n\in\N,\qquad
r_n\to\infty\text{ as }n\to\infty,\\ \forall
x\in\left(\bigcup_{j=1}^\infty \gamma^j\right)\cup\gamma:\quad
d(x_0,x)\ne r_n,\\ k_1:=|\gamma\cap B(x_0,r_1)|>0,\qquad
k_n:=|\gamma\cap (B(x_0,r_n)\setminus B(x_0,r_{n-1}))|>0, \qquad
n\ge2.
\end{gather*}
Since $\gamma_j\to\gamma$ vaguely, we then  conclude that there
exist $j_1\in\N$ such that $$|\gamma^j\cap B(x_0,r_1)|=k_1\qquad
\text{for all } j\ge j_1,$$ and  a numeration of the points of
$\gamma^j\cap B(x_0,r_1)$, $j\ge j_1$, such that $$\gamma^j\cap
B(x_0,r_1)=\{x_k^j\}_{k=1}^{k_1},\qquad x_k^j\to x_k\text{ as
$j\to\infty$ for all }k=1,\dots,k_1.$$ Next, there exist
$j_2\in\N$, $j_2>j_1$, such that
$$|\gamma^j\cap(B(x_0,r_2)\setminus B(x_0,r_1))|=k_2\qquad
\text{for all }j\ge j_2,$$ and  a numeration of the points of
$\gamma^j\cap(B(x_0,r_2)\setminus B(x_0,r_1))$, $j\ge j_2$, such
that
\begin{gather*} \gamma^j\cap(B(x_0,r_2)\setminus
B(x_0,r_1))=\{x_k^j\}_{k=k_1+1}^{k_1+k_2}, \\ x_k^j\to x_k\text{
as $j\to\infty$ for all } k=k_1+1,\dots,k_1+k_2.\end{gather*}
Continuing this procedure by induction, we get for each $j_n\le
j<j_{n+1}$ a numeration $\{x_k^j\}_{k=1}^{k_1+\dots+k_n}$ of the
points of $\gamma^j\cap B(x_0,r_n)$. For such $j$ we choose an
arbitrary numeration $\{x_k^j\}_{k=k_1+\dots+k_n+1}^\infty$ of
$\gamma^j\cap B(x_0,r_n)^{\mathrm c}$, and therefore obtain a
numeration $\{x_k^j\}_{k=1}^\infty$ of $\gamma_j$. Since
$j_n\to\infty$, we thus have a  numeration of all  $\gamma^j$ with
$j\ge j_1$ (for the first $j_1-1$ configurations, we again take
an arbitrary numeration). Now, for any fixed $l\in\N$, take the
minimal $n(l)\in\N$ satisfying $x_l\in B(x_0,r_{n(l)})$ (this
$n(l)$ always exists since $r_n\to\infty$). Then,
$k_1+\dots+k_{n(l)-1}+1\le l\le k_1+\dots+k_{n(l)}$ (where
$k_1+\dots+k_{n(l)-1}:=0$ if $n(l)=1$). By induction, the sequence
$(x_k^j)_{j=j_{n(l)}}^\infty$ converges to $x_k$ as $j\to\infty$
for each $k$ satisfying $k_1+\dots+k_{n(l)-1}+1\le k\le
k_1+\dots+k_{n(l)}$, which immediately yields that the sequence
$(x_l^j)_{j=1}^\infty$ converges to $x_l$ as $j\to\infty$.

According to \eqref{kj}, \eqref{709} is, therefore,  equivalent to
the convergence
\begin{equation}\label{711}
\sum_{(i_1,\dots,i_n)\in\tilde{\N}{}^n}\,\prod_{k=1}^n
p_{t}(x_{i_k}^j,y_k)\to\sum_{(i_1,\dots
,i_n)\in\tilde{\N}{}^n}\prod_{k=1}^n p_{t}(x_{i_k},y_k)\qquad
\text{as }j\to\infty
\end{equation} for each fixed $(y_1,\dots,y_n)\in \tilde X^n$,
$n\in\N$.

We now claim that, for any fixed $\eps>0$, there exist $J,K\in\N$
such that \begin{equation}\label{nuodin} \forall j>J:\quad
\sum_{k=K+1}^\infty\exp\big[-d(x_0,x_k^j)\,\big]<\eps,\qquad
\sum_{k=K+1}^\infty\exp\big[-d(x_0,x_k)\,\big]<\eps.\end{equation}
Indeed, choose any $K\in\N$ such that
\begin{equation}\label{nudva}\sum_{k=K+1}^\infty\exp\big[-d(x_0,x_k)
\,\big]<\frac\eps3,\end{equation} then choose any $J_1\in \N$ such
that
\begin{equation}\label{nutri}\forall j>J_1:\quad
\bigg|\sum_{k=1}^K\exp\big[-d(x_0,x_k^j)\,\big]
-\sum_{k=1}^K\exp\big[-d(x_0,x_k)\,\big]\bigg|<\frac\eps3
,\end{equation} and finally take any $J_2\in\N$ such that
\begin{equation}\label{nuchetyre}\forall j>J_2:\quad \bigg|\sum_{k=1}^\infty
\exp\big[-d(x_0,x_k^j)\,\big]-\sum_{k=1}^\infty\exp\big[-d(x_0,x_k)\,
\big]\bigg|<\frac\eps3.\end{equation} Then, it follows from
\eqref{nudva}--\eqref{nuchetyre} that \eqref{nuodin} holds with
$K$ as above and $J:=\max\{J_1,J_2\}$.

Now, we  conclude from (C1) that, for each fixed $t>0$ and $y\in
X$, there exists $\operatorname{const}_{t,y}>0$ such that
\begin{equation}\label{nupiat}p_t(x,y)\le\operatorname{const}_{t,y}
\exp\big[-d(x_0,x)\,\big].\end{equation} Thus,  from \eqref{710},
\eqref{nuodin}, and \eqref{nupiat}, we easily derive that
\begin{equation}\label{nushest}\sum_{k=1}^\infty \big|p_t(x_k^j,y)- p_t(x_k,y)\big|\to0
 \qquad \text{as
}j\to\infty.\end{equation} Thus, \eqref{711} holds for $n=1$.

Next, we  show that \eqref{711} holds for $n=2$, i.e.,
\begin{equation}\label{nusem} \sum_{i_1=1}^\infty\bigg|
p_t(x_{i_1}^j,y_1)\sum_{i_2\in\N,\, i_2\ne i_1}p_t(x_{i_2}^j,y_2)-
 p_t(x_{i_1},y_1)\sum_{i_2\in\N,\, i_2\ne
i_1}p_t(x_{i_2},y_2)\bigg|\to0\qquad \text{as
}j\to\infty.\end{equation} It follows from \eqref{710} and
\eqref{nushest} that, for each $i_1\in\N$,
\begin{equation}\label{nuvosem}
p_t(x_{i_1}^j,y_1)\sum_{i_2\in\N,\, i_2\ne
i_1}p_t(x_{i_2}^j,y_2)\to p_t(x_{i_1},y_1)\sum_{i_2\in\N,\, i_2\ne
i_1}p_t(x_{i_2},y_2)\qquad \text{as }j\to\infty.\end{equation}
Moreover, we get from  \eqref{nushest} that
\begin{equation}\label{nudeviat}\sum_{i_2\in\N,\, i_2\ne
i_1}p_t(x_{i_2}^j,y_2)\le\sum_{i_2=1}^\infty
p_t(x_{i_2}^j,y_2)\le\operatorname{const}\qquad \forall
j\in\N.\end{equation} Thus, we obtain \eqref{nusem} from
 \eqref{nuodin},  \eqref{nupiat}, \eqref{nuvosem}, and
\eqref{nudeviat}.

Continuing this way, by induction we prove \eqref{711} for each
$n\in\N$.

By virtue of  the majorized convergence theorem, it still remains
to verify that all the functions $G(\cdot)k_t(\gamma^j,\cdot)$
 are majorized by a
function from $L^1(\lambda_m)$. But it follows from
\eqref{esrdrttud} that $$ |k_t(\gamma^j,\{y_1,\dots,y_n\}) |\le
\operatorname{const}_t^n\exp\big[\,d(x_0,y_1)+\dots+d(x_0,y_n)\,\big].$$
 Thus, \eqref{705} implies the
assertion of  the theorem.\quad$\blacksquare$ \vspace{2mm}

We have also the following theorem.

\begin{th}\label{Meller_invariant}
Let $\Bbb D$ denote the set of all measurable functions
$G=(G^{(n)})_{n=0}^\infty$ on $\Gamma_{X,0}$ such that there exist
$\eps=\eps(G)$\rom, $C=C(G)>0$ such that
\begin{gather}\big|G^{(n)}\circ
T_n(x_1,\dots,x_n)\big|\le
C^n\exp\big[-(1+\eps)(d(x_0,x_1)+\dots+d(x_0,x_n))\,\big],\label{klghas}\\
(x_1,\dots,x_n)\in\tilde X^n,\ n\in\N.\notag\end{gather} We define
$${\bf D}:=\big\{\,(KG)\restriction\Gamma_\infty\mid G\in{\Bbb
D}\,\big\}.$$ Then\rom, \begin{equation}\forall p\ge1:\quad {\bf
D}\subset L^p(\Gamma_\infty,\pi_m)\label{awhjrth}\end{equation}
and \begin{equation} \forall t>0:\quad{\bf P}_t{\bf D }\subset{\bf
D}.\label{opewxdki}\end{equation} Furthermore\rom,  for each
$F\in{\bf D}$\rom, ${\bf P}_tF$ is a continuous function on
$\Gamma_\infty$ with respect to the metric $d_1$\rom.\end{th}

\noindent{\it Proof}. By virtue of \eqref{lpseaq}, \eqref{klghas}
implies \eqref{705}, and hence by Theorem~\ref{Meller_Skorohod},
${\bf P}_tF$ is a continuous function  on $\Gamma_\infty$ with
respect to  $d_1$ for each $F\in{\bf D }$.

Next, ${\Bbb D}\subset L^1(\Gamma_{X,0},\lambda_m)$, and
therefore, by Proposition~\ref{prop3.2},
\begin{equation}\label{klxe}{\bf D}\subset
L^1(\Gamma_\infty,\pi_m).\end{equation} For each
$\eta\in\Gamma_{X,0}$, $\sum_{\theta\subset\eta}1= 2^{|\eta|}$,
which yields that
\begin{equation}\label{klcftz}\sum_{(\theta_1,\theta_2,\theta_3)\in{\cal P}_3(\eta)}1\le
6\cdot 4^{|\eta|}.\end{equation} By  \eqref{klcftz} and definition
\eqref{starproduct}, we get that $G_1\star G_2\in{\Bbb D}$ for
arbitrary $G_1,G_2\in{\Bbb D}$. Consequently, by \eqref{waerfed},
\begin{equation}\label{klawht}F_1\cdot F_2\in{\bf D},\qquad
F_1,F_2\in{\bf D}.\end{equation} From \eqref{klxe} and
\eqref{klawht}, we get \eqref{awhjrth}.

Finally, let us show that the set $\bf D$ is invariant under the
action of ${\bf P}_t$. By  \eqref{chyk} and \eqref{awdrawdr}, we
have
\begin{equation}\label{klseawjksre} ({\bf
P} _t(KG))(\gamma)=(K(\tilde{\bf P}_tG))(\gamma),\qquad
\gamma\in\Gamma_\infty,\ G\in{\Bbb D},
\end{equation}
where $\tilde {\bf P}_tG=((\tilde {\bf P}_tG)^{(n)})_{n=0}^\infty$
is given by
\begin{gather}(\tilde {\bf P}_t G)^{(0)}= G^{(0)},\notag\\
(\tilde {\bf P}_tG) ^{(n)}\circ T_n(x_1,\dots,x_n)=\int_{\tilde
X^n}G^{(n)}\circ T_n(y_1,\dots,y_n)\,p_{t,x_1}(dy_1)\dotsm
p_{t,x_n}(dy_n),\notag\\ (x_1,\dots,x_n)\in\tilde X^n,\
n\in\N.\notag\end{gather}

Now, by using (C1), we derive, for  $G=(G^{(n)})_{n=0}^\infty$
satisfying \eqref{klghas} with some $\eps>0$ and $C>0$,
\begin{gather*}\big|(\widetilde
\P_tG)^{(n)}\circ T_n(x_1,\dots,x_n)\big|\le \\ \le \int_{X^n}C^n
\exp\big[-(1+\eps)(d(x_0,y_1)+\dots+d(x_0,y_n))\,\big]\\\times
(C_t')^n\exp\big[-(1+\myfrac\eps2)(d(x_1,y_1)+\dots+d(x_n,y_n))\,\big]
\,m(dy_1)\dotsm m(dy_n)\\\le (CC_t')^n
\exp\big[-(1+\myfrac\eps2)(d(x_0,x_1)+\dots+d(x_0,x_n))\,\big]
\bigg(\int_X\exp\big[-
\myfrac\eps2\,d(x_0,y)\,\big]\,m(dy)\bigg)^n.\end{gather*} Hence,
because of  \eqref{lpseaq}, $\tilde {\bf P}_tG\in{\Bbb D}$, and
\eqref{opewxdki} follows from \eqref{klseawjksre}.\quad
$\blacksquare$\vspace{2mm}

For an arbitrary measurable, bounded, symmetric function
$G^{(n)}(x_1,\dots ,x_n)$ on $X^n$, $n\in\N$, with  bounded
support, one can introduce the following monomial on
$\Gamma_\infty$ of  $n$-th order with  kernel $G^{(n)}$:
$$\Gamma_\infty\ni\gamma\mapsto\la G^{(n)},{:}\,\gamma^{\otimes
n}\,{:}\ra:= \sum_{\{x_1,\dots,x_n\}\subset\gamma}
G^{(n)}(x_1,\dots,x_n)=(KG^{(n)})(\gamma).$$ It is natural to call
a finite sum of functions of such type and a constant a
cylindrical polynomial on $\Gamma_\infty$ with bounded
coefficients. We denote by ${\cal FP}_{\mathrm bc}(\Gamma_\infty)$
the set of all such polynomials on $\Gamma_\infty$. Thus, ${\cal
FP}_{\mathrm bc}(\Gamma_X)$ is nothing but the image of the set
$B_{\mathrm bs}(\Gamman)$ under the $K$-transform.

 Since every function $G\in B_{\mathrm
bs}(\Gamman)$ satisfies \eqref{klghas}, we get the following
consequence of Theorem~\ref{Meller_invariant}:

\begin{cor}\label{klaseshj}
We have the inclusion ${\cal FP}_ {\mathrm
bc}(\Gamma_{\infty})\subset {\bf D} $\rom. In particular\rom, for
each polynomial $F\in {\cal FP}_ {\mathrm
bc}(\Gamma_{\infty})$\rom, the function ${\bf P}_tF$ is continuous
on $\Gamma_{\infty}$ with respect to the metric $d_{1}$\rom.
\end{cor}

A measurable function $F\colon\Gamma_\infty\to\R$ is called local
if there exist an open, relatively  compact set $\Lambda\subset X$
and a measurable function $\tilde F:\Gamma_{\Lambda,0}\to\R$ such
that $F(\gamma)=\tilde F(\gamma_\Lambda)$ for all
$\gamma\in\Gamma_\infty$. The following corollary is an analog of
the classical strong Feller property for the heat semigroup on the
configuration space.

\begin{cor}\label{awnmtezd} Each measurable bounded local
function $F\colon\Gamma_\infty\to\R$ belongs to $\bf D$\rom, and
hence the function ${\bf P}_{t}F$ is continuous on $\Gamma_\infty$
with respect to the metric $d_1$\rom.\end{cor}

\noindent{\it Proof}. Let $F(\gamma)=\tilde F(\gamma_\Lambda)$
with $\Lambda$ and $\tilde F$ as above. By
\cite[Proposition~3.5]{Kuna1}, one can  explicitly calculate the
inverse $K$-transform of $F$: \begin{equation}\label{werrf}
(K^{-1}F)(\eta)=\begin{cases}\sum\limits_{\theta\subset\eta}(-1)^{|\eta\setminus\theta|}
\,\tilde F(\theta),&\text{if }\eta\in\Gamma_{\Lambda,0},\\
0,&\text{otherwise}.\end{cases}\end{equation} Set
$G=(G^{(n)})_{n=0}^\infty:=K^{-1}F$. Let $C:=\sup|F|$. Since
$\sum_{\theta\subset\eta}1=2^{|\eta|}$, we conclude from
\eqref{werrf} that, for each $n\in\N$, $\{G^{(n)}\ne
0\}\subset\Gamma_{\Lambda}^ {(n)}$ and $|G^{(n)}|$ is bounded by
the constant $C2^n$. Therefore, $G\in{\Bbb D}$ and $F=KG\in{\bf D
}$.\quad$\blacksquare$

\begin{rem}\rom{Let $t>0$. As a consequence of
Corollary~\ref{awnmtezd} we have that any Markovian kernel $\tilde
{\bf P}_t$ on $(\Gamma_\infty,{\cal B}(\Gamma_\infty))$ such that
$\tilde {\bf P}_tF$ is continuous with respect to $d_1$ and
$\tilde {\bf P}_tF$ is a $\pi_m$-version of $e^{-tH^\Gamma}F$ for
each measurable bounded local function $F:\Gamma_\infty\to\R$ must
coincide  with ${\bf P}_t$. This follows from the fact that
$\pi_m(U)>0$ for any nonempty  set $U\in{\cal B}(\Gamma_\infty)$
which is open in the topology generated by the metric $d_1$. The
latter can be proved as follows. For any fixed $\hat
\gamma\in\Gamma_\infty$ and $\eps>0$, there exist $R>0$ and
$\delta>0$ such that for each $r\ge R$
\begin{align}&\pi_m\big(
d_1(\gamma,\hat\gamma)<\eps\big)\notag\\ &\qquad \ge
\pi_m\big(d_{\mathrm
K}(\gamma_{\Lambda_r},\hat\gamma_{\Lambda_r})<\delta,\
|B_1(\gamma_{\Lambda_r})-B_1(\hat\gamma_{\Lambda_r})|<\delta,\
B_1(\gamma_{\Lambda_r^{\mathrm c}})<\delta\big)\notag \\
&\qquad=\pi_m\big(d_{\mathrm
K}(\gamma_{\Lambda_r},\hat\gamma_{\Lambda_r})<\delta,\
|B_1(\gamma_{\Lambda_r})-B_1(\hat\gamma_{\Lambda_r})|<\delta\big)\,
\pi_m\big(B_1(\gamma_{\Lambda_r^{\mathrm
c}})<\delta\big).\label{sezgu}\end{align} Here, $\Lambda_r:=\{x\in
X:d(x_0,x)<r\}$ and the functional $B_1$ is defined on the space
$\Gamma_{X,0}$ by  the same formula \eqref{judr}.
 The first factor in
\eqref{sezgu} is obviously positive, while the positivity of the
second factor for sufficiently big $r>0$ is implied by
$$1=\pi_m\big(B_1(\gamma\big)<\infty)=\pi_m\bigg(\bigcup_{r=1}^\infty\{
B_1(\gamma_{\Lambda_r^{\mathrm c}})<\delta\}
\bigg)=\lim_{r\to\infty} \pi_m\big(B_1(\gamma_{\Lambda_r^{\mathrm
c}})<\delta\big).$$ }\end{rem}

\section{Feller property of the heat semigroup
with respect to the intrinsic metric of the Dirichlet form }

For presenting another type of Feller property of the heat
semigroup $(e^{-tH^\Gamma})_{t\in\R_+}$, we shall need the space
$\dd\Gamma_X$ of all $\Z_+$-valued Radon measures $\gamma$ on $X$
such that $\gamma(X)=\infty$. This space is the closure of
$\Gamma_X$ in the $d_{\mathrm K}$ metric.  The space $\dd\Gamma_X$
is equipped with the topology induced by the vague topology on
${\cal M}(X)$, and let ${\cal B}(\dd\Gamma_X)$ denote the
corresponding Borel $\sigma$-algebra.

Furthermore, let $$\overset{{.}{.}}{\Gamma}_n:=\big\{\,\gamma\in
\overset{{.}{.}}{\Gamma}_X: B_n(\gamma)<\infty\,\big\},$$ where
$B_n$ is as in \eqref{judr}, but defined on all of $\dd\Gamma_X$.
Let $\overset{{.}{.}}{\Gamma}_\infty:=\bigcap_{n=1}^\infty
\overset{{.}{.}}{\Gamma}_n$. We extend the metric $d_\infty$ to
$\dd\Gamma_\infty$ using the same formula \eqref{secondmetric}.
The Borel $\sigma$-algebra ${\cal B}(\dd\Gamma_\infty)$
corresponding to the $d_\infty$ metric  coincides with the trace
$\sigma$-algebra of ${\cal B}(\dd\Gamma_X)$ on $\dd\Gamma_\infty$.

Let $\widehat X^\N$ denote the (${\cal B}(X^\N)$-measurable)
subset of $X^\N$ consisting of those $(x_1,x_2,\dots)\in X^\N$ for
which the number of the $x_k$'s in any compactum in $X$ is finite.
Evidently, one can identify $\dd\Gamma_X$ with the factor space $
\widehat X^\N/S_\infty$. Analogously to \eqref{oizgdr}, we define
the corresponding quotient map  $I\colon \widehat X^\N\to
\overset{{.}{.}}{\Gamma}_X$ by
\begin{equation}\label{awghzj}\widehat X^\N\ni{\bf x}=(x_1,x_2,\dots)\to
I{\bf x}:=[x_1,x_,\dots]\in
\overset{{.}{.}}{\Gamma}_X,\end{equation}  which is measurable, as
can be seen by similar arguments as those following
\eqref{oizgdr}.

For each $t>0$ and $\gamma\in\dd\Gamma_X$, we define the measure
$\P_{t,\gamma}$ on $\dd\Gamma_X$ as the image under the mapping
\eqref{awghzj} of the restriction to $\widehat X^\N$ of any
measure ${\Bbb P}_{t,{\bf x}}:=\bigotimes_{k=1}^\infty p_{t,x_k}$,
${\bf x }=(x_k)_{k=1}^\infty\in I^{-1}\{\gamma\}$  (the resulting
measure is independent of the choice of ${\bf x }\in
I^{-1}\{\gamma\}$). Thus, by Lemma~\ref{lem1}, $\P_{t,\gamma}$ is
either   a probability measure  or
 zero measure on $\dd\Gamma_X$ depending on  whether the series
  $\sum_{k=1}^\infty p_{t,x_k}(\Lambda)$ converges for each
  compact $\Lambda\subset X$, or not, and $\P_{t,\gamma}(\dd\Gamma_X)=\P_{t,\gamma}
  (\Gamma_X)$. In the same way as we proved
  Proposition~\ref{quadrofonia}, we conclude that, for each
  $\gamma\in\dd\Gamma_\infty$,
  $$\P_{t,\gamma}(\dd\Gamma_\infty)=\P_{t,\gamma}(\Gamma_\infty)=1.$$

Following \cite{Ro}, we introduce the  $L^2$-Wasserstein type
distance $\rho$ on $\dd\Gamma_X$ setting, for any
$\gamma_1=[x_1,x_2,\dots ]$ and $\gamma_2=[y_1,y_2,\dots]$ from
$\dd\Gamma_X$, \begin{equation}\label{def}
\rho(\gamma_1,\gamma_2):=\inf\bigg\{ \bigg( \sum_{k=1}^\infty
d(x_k,y_{\sigma(k)})^2 \bigg)^{1/2}\,{\Big|}\,\sigma\in S_\infty
\bigg\}.
\end{equation}
 Notice that $\rho$ is a
pseudo-metric, i.e., it takes values in $[0,\infty]$. Obviously,
convergence with respect to $\rho$ implies vague convergence. We
recall that $\rho$ is the intrinsic metric of the Dirichlet form
obtained as the closure of \eqref{three}, see \cite{Ro}.

Analogously to \eqref{fuka}, we set for a  measurable function $F$
on $\dd\Gamma_X$:
\begin{equation}\label{cheese}
(\Pe_t
F)(\gamma):=\int_{\dd\Gamma_X}F(\xi)\,\Pe_{t,\gamma}(d\xi),\qquad
\gamma\in\dd\Gamma_X,\, t>0,
\end{equation}
provided the integral on the right hand side of \eqref{cheese}
exists.

In the rest of this section, we shall be concerned only with the
case $X=\R^d$. Let us recall that the heat kernel has now the form
\begin{equation}\label{heatkernel}p(t,x,y)=(4\pi t)^{-d/2}\big[-\myfrac1{4t}\,
|x-y|^2\,\big].\end{equation} We shall show that the space $\Crb$
of all bounded functions on $\dddGamma$ which are continuous with
respect to the  $\rho$ metric is invariant under ${\bf P}_t$ for
all $t>0$.

\begin{th} \label{koryto}
We have\rom: $$\Pe_t(\Crb)\subset \Crb,\qquad t>0.$$
\end{th}

\noindent{\it Proof}. First, we note that the distance $\rho$ can
be extended from $\dddGamma=\widehat{\R^d}{} ^\N/S_\infty$  to the
bigger space $(\R^d)^\N/S_\infty$ by using the same formula
\eqref{def} for calculating the distance between any
$\gamma_1=[x_1,x_2,\dots]$ and $\gamma_2=[y_1,y_2,\dots]$ from
$(\R^d)^\N/S_\infty$.

It follows directly from \eqref{def} that, if $\gamma_1$ and
$\gamma_2$ are two elements of $(\R^d)^\N/S_\infty$ having finite
$\rho$ distance, then   $\gamma_1\in\dddGamma$ implies
$\gamma_2\in \dddGamma$. Therefore, any function $F\in \Crb$ can
be extended to a continuous bounded function on
$(\R^d)^\N/S_\infty$, again denoted by $F$, as follows:
$$F\restriction ((\R^d)^\N/S_\infty)\setminus\dddGamma:=0.$$ Then,
\eqref{cheese} yields
\begin{equation}\label{chren}
(\Pe_tF)(\gamma)=\int_{(\R^d)^\N}
F(y_1,y_2,\dots)\,\bigotimes_{k=1}^\infty p_{t,x_k}(dy_k),\qquad
\gamma=[x_1,x_2,\dots]\in\dddGamma.\end{equation} Notice that, in
this formula, $F$ is considered as an $S_\infty$-invariant
function on $(\R^d)^\N$.

Since $F$ is bounded, so is the function $\Pe_tF$, and we only
have to prove the continuity. To this end, let $\gamma^j\to
\gamma$ in the $\rho$ metric.
 By Lemma~4.1 in \cite{Ro}, there always
exists a representative $(x_k^j)_{k=1} ^\infty$ of $\gamma^j$ such
that
\begin{equation}\label{afrika}
\rho(\gamma^j,\gamma)=\bigg(\sum_{k=1}^\infty
|x_k^j-x_k|^2\bigg)^{1/2},\qquad j\in\N.\end{equation}
\eqref{heatkernel} and \eqref{chren}   imply
\begin{gather}
(\Pe_t F)(\gamma^j)=\int_{(\R^d)^\N}
F(y_1+x_1^j,y_2+x_2^j,\dots)\,\bigotimes_{k=1}^\infty
p_t(dy_k),\label{redka}\\ p_t(dy)=(4\pi
t)^{-d/2}\exp\big[-\myfrac1{4t}\,|y|^2\,\big]\,dy.\notag\end{gather}
Since the integrand in \eqref{redka} is a bounded function, it
suffices to show that for any fixed $(y_k)_{k=1}^\infty\in
(\R^d)^\N$ $$F(y_1+x_1^j,y_2+x_2^j,\dots)\to
F(y_1+x_1,y_2+x_2,\dots)\quad\text{as }j\to\infty.$$ But this
follows from the fact that $F$ is continuous in the $\rho$ metric
and from the convergence $$ \rho\big(
[y_1+x_1^j,y_2+x_2^j,\dots],[y_1+x_1,y_2+x_2,\dots]
\big)\le\bigg(\sum_{k=1}^\infty|y_k+x_k^j-y_k-x_k|^2\bigg)^{1/2}
=\rho(\gamma_j,\gamma)\to0$$ as $j\to\infty$, which, in turn, is
implied by \eqref{afrika}.\quad $\blacksquare$

\begin{rem}\rom{Theorem~\ref{koryto}, in particular,  yields that if
$\dddGamma$ is of full $\Pe_{t,\gamma}$ measure for some
$\gamma\in\dddGamma$, then it is also of full $\Pe_{t,\gamma'}$
measure for each $\gamma'\in\dddGamma$ such that $
\rho(\gamma,\gamma')<\infty$. For it suffices to note that $
1\in\Crb$ and $(\Pe_t 1)(\gamma)=\Pe_{t,\gamma} (\ddGamma)$.}
\end{rem}

Finally, we shall present another version of the latter  theorem.
Let $C_{\rho,{\mathrm b}}(\overset{{.}{.}}{\Gamma}_{\infty})$
denote the space of all bounded functions on
$\overset{{.}{.}}{\Gamma}_{\infty}$ that are continuous with
respect to the $\rho$ metric.  It is easy to see that
 each element of $(\R^d)^ \N/S_\infty$ having a
finite distance to $\overset{{.}{.}}{\Gamma}_{\infty}$ itself
belongs to $\overset{{.}{.}}{\Gamma}_{\infty}$. Therefore, any
function  $F\in C_{\rho,{\mathrm
b}}(\overset{{.}{.}}{\Gamma}_{\infty})$ can be extended to a
function from $\Crb$ by setting $F$ to be equal to zero on
$\dddGamma \setminus \overset{{.}{.}}{\Gamma}_{\infty}$. We note
that the convergence on $\overset{{.}{.}}{\Gamma}_{\infty}$ with
respect to the $\rho$ metric implies the convergence with respect
to the $d_\infty$ metric.

Since for each   $\gamma\in \overset{{.}{.}}{\Gamma}_{\infty}$ the
measure ${\bf P}_{t,\gamma}$ is concentrated on
$\overset{{.}{.}}{\Gamma}_{\infty}$, we get the following
corollary of Theorem~\ref{koryto}:

\begin{cor}\label{styk}
We have\rom: $$\Pe_{t}(C_{\rho,{\mathrm
b}}(\overset{{.}{.}}{\Gamma}_{\infty}))\subset C_{\rho,{\mathrm
b}}(\overset{{.}{.}}{\Gamma}_{\infty}),\qquad t>0.$$
\end{cor}

\section{ Brownian motion on the configuration space}
We again consider the case of a general manifold $X$. Analogously
to \eqref{label1}, \eqref{label2}, we define the family of kernels
$(\PP_t)_{t\in\R_+}$ on the space
$(\overset{{.}{.}}\Gamma_\infty,{\cal
B}(\overset{{.}{.}}\Gamma_\infty))$  setting
$$\PP_t(\gamma,A):=\P_{t,\gamma}(A),\qquad
\gamma\in\overset{{.}{.}}\Gamma_\infty,\,A\in{\cal
B}(\overset{{.}{.}}\Gamma_\infty),\,t\in\R_+,$$ where ${\bf
P}_{t,\gamma}$, $t>0$, $\gamma\in\dd\Gamma_X$, is defined as in
the previous section and  $$\P_{0,\gamma}:=\varepsilon_\gamma.$$
Analogously to Proposition~\ref{wadr}, we conclude that
$(\P_t)_{t\in\R_+}$ is
 a Markov semigroup of kernels on
$(\overset{{.}{.}}\Gamma_\infty,{\cal
B}(\overset{{.}{.}}\Gamma_\infty))$.

In this section, we shall give a direct construction of the
independent infinite particle process. Under some additional
conditions on the manifold $X$, we shall show that the resulting
process  is the unique continuous Markov process on
$\dd\Gamma_\infty$ with transition probabilities
$\P_{t}(\gamma,\cdot)$. (We note that we are forced to deal with
the space $\dd\Gamma_\infty$, rather than $\Gamma_\infty$, because
in the general case we cannot exclude collision of the particles,
see Corollary~\ref{diagonal} below).

First, we strengthen a little bit condition (C1) by requiring the
following stronger upper bound:

\begin{description}
\item[\rom{(C$1'$)}] For each $t>0$, there exist
$\vartheta_t\in(0,t)$, $C_t>0$, and  $\eps_t>0$ such that
$$p(s,x,y)\le C_t\exp\big[-d(x,y)^{1+\eps_t}\,\big], \qquad
s\in(t-\vartheta_t,t+\vartheta_t),\ x,y\in X.$$ \end{description}

Evidently, \eqref{seghdr} and \eqref{waeergdzt} imply (C$1'$).

Let us introduce the function
$$\tau(\delta,r):=\sup_{t\in(0,\delta]}\,\sup_{x\in
X}\int_{B(x,r)^{\mathrm c }}p(t,x,y)\,m(dy),\qquad \delta>0,\,
r>0.$$ Because of \eqref{6745635}, $\tau(\delta,r)\le1$ for all
$\delta,r>0$, and for each fixed $r>0$ $\tau(\cdot,r)$ is an
increasing function on $(0,\infty)$.

 Let
$\Omega:=C(\R_+;X)$ denote the space of all continuous functions
(paths) from $\R_+$ to $X$, and let $\cal F$ be the product
$\sigma$-algebra on $\Omega$, i.e.,
\begin{equation}\label{347845}{\cal
F}:=\sigma\{x_t,\,t\in\R_+\},\end{equation} where
$\Omega\ni\omega\mapsto x_t(\omega):=\omega(t)\in X$. For each
$x\in X$, let  $P_x$ denote the
 measure on $(\Omega,{\cal F})$ corresponding to  Brownian
motion on $X$ starting at $x$.

We shall need the following lemma.

\begin{lem}\label{alaNelson} Let $0\le a<b$ with $b-a\le \delta$\rom.
Then\rom, for each $x\in X$ and $r>0$\rom, $$P_x\big(\exists
s,t\in[a,b]:\,
d(\omega(s),\omega(t))>r)\big)\le2\tau(\delta,\myfrac14\,r).$$
\end{lem}

\noindent{\it Proof}. This lemma is a straightforward
generalization of \cite[Appendix~A, Lemma~4]{Nelson}, which deals
with the usual Brownian motion on $\R^d$. However, for
completeness, we present a proof of this lemma in the Appendix.
\quad $\blacksquare$ \vspace{2mm}

We  suppose:

\begin{description}\item[\rom{(C3)}] For each fixed $r>0$,
\begin{equation}\label{C31}\tau(\delta,r)\to 0\qquad \text{as }\delta\to0,\end{equation}
 and there
exist $\tilde\delta>0$ and $C>0$ such that
\begin{equation}\label{C32}\tau(\tilde\delta,r)\le C e^{-r},\qquad
r>0.\end{equation}
\end{description}

The following simple lemma gives a sufficient condition for (C3)
to hold.

\begin{lem}\label{drters} Suppose that the manifold $X$ has nonnegative
Ricci curvature and  the heat kernel $p(t,x,y)$ of $X$ satisfies
the Gaussian upper bound for small values of $t$\rom:
\begin{equation}\label{esreserew}p(t,x,y)\le
Ct^{-n/2}\,\exp\bigg[-\frac{d(x,y)^2}{Dt}\,\bigg],\qquad
t\in(0,\tilde\delta],\; x,y\in X,\end{equation} where $n\in\N$ and
$\tilde\delta$\rom, $C$\rom, and $D$ are positive constants\rom.
Then\rom, \rom{(C3)} is satisfied\rom.
\end{lem}

\begin{rem}\label{rtfse}\rom{Concerning the Gaussian upper bound
\eqref{esreserew}, see e.g.\ \cite{Davies,Gr1,Gr2} and the
references therein. In particular, \eqref{esreserew} is implied by
the estimate $$p(x,x,t)\le C t^{-n/2},\qquad t>0,\; x\in X.$$
}\end{rem}

\noindent{\it Proof of Lemma\/} \ref{drters}. Fix $r>0$, then for
$t\in(0,\tilde\delta]$ and $x\in X$ we get, by \eqref{scuople} and
\eqref{esreserew},
\begin{gather} \int_{B(x,r)^{\mathrm c}} p(t,x,y)\,m(dy)\le\int_{B(x,r)^{\mathrm c}}
C t^{-n/2}\,\exp\bigg[-\frac{d(x,y)^2}{Dt}\,\bigg]\, m(dy)\notag\\
\le\operatorname{const}_1\int_{B(x,r)^{\mathrm c}}
\exp\bigg[\frac1t\bigg(\frac{r^2}{2D}-\frac{d(x,y)^2}{D}\bigg)\bigg]\,m(dy)\notag\\
\le\operatorname{const}_1 \sum_{n=1}^\infty
\exp\bigg[\frac1t\bigg(\frac{r^2}{2D}-\frac{(r+n-1)^2}{D}\bigg)\bigg]\,
m\big(B(x,r+n)\setminus B(x,r+n-1)\big)\notag\\
\le\operatorname{const}_2 \sum_{n=1}^\infty
\exp\bigg[\frac1t\bigg(\frac{r^2}{2D}-\frac{(r+n-1)^2}{D}\bigg)\bigg]\,(r+n)^d.
\label{jkjkjw}\end{gather} Since each term of the latter series
monotonically converges to zero as $t\to0$, we get \eqref{C31}.

Next, because $\tau(\delta,r)$ is bounded by 1, it is enough to
verify that \eqref{C32} holds for all $r\ge R$ with some $R>0$.
Now, analogously to \eqref{jkjkjw}, we get for each
$t\in(0,\tilde\delta]$ and $r\ge1$
\begin{gather*} \int_{B(x,r)^{\mathrm
c}}p(t,x,y)\,m(dy)\le\operatorname{const}_3\sum_{n=1}^\infty\exp\bigg[
\frac1{\tilde\delta}\,\bigg(\frac1{2D}-\frac{(r+n-1)^2}D\bigg)\bigg](r+n)^d\\
\le\operatorname{const}_4\sum_{n=1}^\infty
\exp\big[-2(r+n-1)+(r+n)\big] = \operatorname{const}_4\,
e^{-r}\sum_{n=1}^\infty e^{-n},\end{gather*} which yields the
statement.\quad$\blacksquare$

\begin{th}\label{klhjuzip} Let \rom{(C$1'$), (C2),} and \rom{(C3)} hold\rom.
Then\rom, the independent infinite particle process can be
realized as the unique  continuous\rom, time homogeneous Markov
process $${\bf M}=({\pmb{ \Omega}},{\bf F},({\bf
F}_t)_{t\in\R_+},({\pmb \theta}_t)_{t\in\R_+}, ({\bf P
}_\gamma)_{\gamma\in\dd\Gamma_\infty},({\bf X}_t)_{t\in\R_+})$$ on
the state-space $(\dd\Gamma_\infty,{\cal B}(\dd\Gamma_\infty))$
with transition probability function $(\PP_t)_{t\in\R_+}$
\rom(cf\rom.\ e\rom.g\rom{.\ \cite{BG}).}\end{th}

\noindent {\it Proof of Theorem\/} \ref{klhjuzip}. Let us consider
the set  $\Omega^\N$ and the product $\sigma$-algebra ${\cal
C}_\sigma(\Omega^\N)$ on it that is constructed from the
$\sigma$-algebra $\cal F$ on $\Omega$.

We fix any  ${\bf x}=(x_k)_{k=1}^\infty\in\widehat X^\N$ such that
\begin{equation}\label{klaw}\sum_{k=1}^\infty
\exp\big[-\mbox{$\frac1{n}$}\,d(x_0,x_k)\,\big]<\infty\qquad\forall
n\in\N \end{equation} and define the product measure
\begin{equation}\label{oiwegf}{\Bbb P}_{\bf x}:=\bigotimes _{k=1}^\infty P_{x_k}
\end{equation} on
$(\Omega^\N,{\cal C}_\sigma(\Omega^\N))$.

By using \eqref{iikkkf}, we conclude from \eqref{klaw} that
\begin{align*}&\int_{X^\N}\sum_{k=1}^\infty
\exp\big[-\mbox{$\frac1{n}$}\,d(x_0,y_k)\,\big]\,\bigotimes
_{k=1}^\infty p_{t,x_k}(dy_k)\\&\qquad =\sum_{k=1}^\infty \int_X
\exp\big[-\mbox{$\frac1{n}$}\,d(x_0,y)\,\big]\,p_{t,x_k}(dy)<\infty\end{align*}
for all $t>0$ and $n\in\N$, which yields that, for each fixed
$t\in\R_+$,
\begin{align}&\sum\limits_{k=1}^\infty
\exp\big[-\mbox{$\frac1{n}$}\,d(x_0,\omega_k(t))\,\big]<\infty\qquad
\text{for all $n\in\N$ and ${\Bbb P }_{\bf x}$-a.e.\
$\omega=(\omega_k)_{k=1}^\infty\in\Omega^\N$,}\notag\\ &\text{in
particular, $(\omega_k(t))_{k=1}^\infty\in\widehat X{}^\N$ for
such $\omega\in\Omega^\N$.}\label{awuk}\end{align}

\begin{lem}\label{lemmaeins} For each fixed $t\in\R_+$ and
${\bf x}=(x_k)_{k=1}^\infty\in\widehat X^\N$ satisfying
\eqref{klaw}, we have
\begin{equation}\label{yerdztf}{\Bbb P}_{\bf
x}\bigg(\bigcup_{i=1}^\infty \bigcap_{k=1}^\infty\big\{\,
d(\omega_k(s),\omega_k(t))\le\max\{1,\myfrac12\,d(x_0,\omega_k(t))\}\
\forall
s\in(t,t+\myfrac1i]\,\big\}\bigg)=1.\end{equation}\end{lem}

\noindent{\it Proof}.  First, we will prove  \eqref{yerdztf} for
$t=0$. Thus, we have to show that
\begin{equation}\label{yeklgzrdztf}{\Bbb P}_{\bf
x}\bigg(\bigcup_{i=1}^\infty \bigcap_{k=1}^\infty\big\{\,
d(\omega_k(s),x_k)\le\max\{1,\myfrac12\,d(x_0,x_k)\}\ \forall
s\in(0,\myfrac1i]\,\big\}\bigg)=1,\end{equation} or equivalently
\begin{equation}\label{wye}{\Bbb P}_{\bf
x}\bigg(\bigcap_{i=1}^\infty \bigcup _{k=1}^\infty\big\{\,\exists
s\in(0,\myfrac1i]:
d(\omega_k(s),x_k)>\max\{1,\myfrac12\,d(x_0,x_k)\}\,\big\}\bigg)=0.\end{equation}
Since \begin{align}&{\Bbb P}_{\bf x}\bigg(\bigcap_{i=1}^\infty
\bigcup _{k=1}^\infty\big\{\,\exists s\in(0,\myfrac1i]:
d(\omega_k(s),x_k)>\max\{1,\myfrac12\,d(x_0,x_k)\}\,\big\}\bigg)\notag\\
&\qquad=\lim_{i\to\infty}{\Bbb P}_{\bf x}\bigg( \bigcup
_{k=1}^\infty\big\{\,\exists s\in(0,\myfrac1i]:
d(\omega_k(s),x_k)>\max\{1,\myfrac12\,d(x_0,x_k)\}\,\big\}\bigg),\label{rdil}\end{align}
we have, by Lemma~\ref{alaNelson}, \begin{align} &{\Bbb P}_{\bf
x}\bigg( \bigcup _{k=1}^\infty\big\{\,\exists s\in(0,\myfrac1i]:
d(\omega_k(s),x_k)>\max\{1,\myfrac12\,d(x_0,x_k)\}\,
\big\}\bigg)\notag\\&\qquad\le\sum_{k=1}^\infty {\Bbb P}_{\bf
x}\big(\exists
s\in(0,\myfrac1i]:d(\omega_k(s),x_k)>\max\{1,\myfrac12\,d(x_0,x_k)\}\big)
\notag\\&\qquad=\sum_{k=1}^\infty
 P_{ x_k}\big(\exists s\in(0,\myfrac1i]:
d(\omega(s),x_k)>\max\{1,\myfrac12\,d(x_0,x_k\}\big)\notag\\
&\qquad \le2\sum_{k=1}^\infty
\tau\big(\myfrac1i,\max\{\myfrac14,\myfrac18\, d(x_0,x_k)
\}\big).\label{yxmn}\end{align} By  \eqref{C31},
\begin{equation}\label{ljkseghrtdf}
\tau\big(\myfrac1i,\max\{\myfrac14,\myfrac18\, d(x_0,x_k)
\}\big)\to0\qquad \text{as }i\to\infty\end{equation} for each
$x_k$. On the other hand, it follows from \eqref{C32} that, for
any $i\in\N$ satisfying $\frac1i\le\tilde\delta$, the latter
series in \eqref{yxmn} is majorized by the series
$$2C\sum_{k=1}^\infty
\exp\big[-\max\{\myfrac14,\myfrac18\,d(x_0,x_k)\}\,\big],$$ which
converges due to
 \eqref{klaw}. Hence, \eqref{wye}
follows from \eqref{rdil}--\eqref{ljkseghrtdf}, and the monotone
convergence theorem.

Next, using the  Markov property of Brownian motion on $X$, we
easily conclude that
\begin{align}&{\Bbb P}_{\bf x}\bigg(\bigcup_{i=1}^\infty
\bigcap_{k=1}^\infty
\big\{\,d(\omega_k(s),\omega_k(t))\le\max\{1,\myfrac12\,
d(x_0,\omega_k(t))\}\ \forall
s\in(t,t+\myfrac1i]\,\big\}\bigg)\notag\\ &\qquad=\int_{\widehat
X^\N}{\Bbb P}_{t,{\bf x}}(d{\bf y})\,{\Bbb P}_{\bf
y}\bigg(\bigcup_{i=1}^\infty\bigcap_{k=1}^\infty\big\{\,
d(\omega_k(s),y_k)\le\max\{1,\myfrac12\,d(x_0,y_k)\}\ \forall
s\in(0,\myfrac1i]\,\big\}\bigg),\label{klawbn}\end{align} where
${\Bbb P}_{t,{\bf x}}$ is the distribution of
$\omega(t)=(\omega_k(t))_{k=1}^\infty$ under ${\Bbb P}_{\bf x}$.
Now, \eqref{yerdztf} follows from \eqref{awuk},
\eqref{yeklgzrdztf}, and
\eqref{klawbn}.\quad$\blacksquare$\vspace{2mm}

\begin{lem}\label{newlemma} For each fixed $t>0$ and
 ${\bf x}=(x_k)_{k=1}^\infty\in\widehat X^\N$ satisfying
 \eqref{klaw}\rom, we have \begin{equation*}{\Bbb
 P}_{\bf x}\bigg(\bigcup_{i={\cal I}_t}^\infty\bigcap_{k=1}^\infty\big\{\,
 d(\omega_k(t-\myfrac1i),\omega_k(s))\le\max\{1,\myfrac12\,
 d(x_0,\omega_k(t-\myfrac1i))\}\ \forall
 s\in(t-\myfrac1i,t]\,\big\}\bigg)=1,\end{equation*}
 where ${\cal I}_t:=[t^{-1}]+1$ \rom($[a]$ denoting
 the integer part of $a>0$\rom)\rom.\end{lem}

\noindent {\it Proof}. It is enough to show that
\begin{equation}\label{huizgurfrztrft} \lim_{i\to\infty}{\Bbb
P}_{\bf x}\bigg(\bigcup_{k=1}^\infty\big\{\, \exists
s\in(t-\myfrac1i,t]:
d(\omega_k(t-\myfrac1i),\omega_k(s))>\max\{1,\myfrac12\,d(x_0,\omega_k(t-\myfrac1i))\}
\,\big\}\bigg)=0.\end{equation} Using Lemma~\ref{alaNelson}, we
get
\begin{align}&{\Bbb P}_{\bf x}\bigg(\bigcup_{k=1}^\infty\big\{\,
\exists s\in(t-\myfrac1i,t]:
d(\omega_k(t-\myfrac1i),\omega_k(s))>\max\{1,\myfrac12\,d(x_0,\omega_k(t-\myfrac1i))\}
\,\big\}\bigg)\notag\\ &\qquad\le \sum_{k=1}^\infty P_{x_k}\big(
\exists s\in(t-\myfrac1i,t]:
d(\omega(t-\myfrac1i),\omega(s))>\max\{1,\myfrac12\,
d(x_0,\omega(t-\myfrac1i))\}\big)\notag\\&\qquad=\sum_{k=1}^\infty
\int_X m(dy)\, p(t-\myfrac1i,x_k,y)P_y\big(\exists
s\in(0,\myfrac1i]:d(y,\omega(s))>\max\{1,\myfrac12\,d(x_0,y)\}\big)\notag\\
&\qquad\le\sum_{k=1}^\infty \int_X m(dy)\,
p(t-\myfrac1i,x_k,y)\,2\tau(\myfrac1i,\max\{\myfrac14,\myfrac18\,d(x_0,y)\}).
\label{sersre}\end{align}

 Now, it follows from \eqref{C31} that
\begin{equation}\label{12365}
p(t-\myfrac1i,x_k,y)\,2\tau(\myfrac1i,\max\{\myfrac14,\myfrac18\,d(x_0,y)\})
\to0\qquad\text{as }i\to\infty\end{equation} for each  $x_k$ and
$y\in X$. Next, by (C$1'$) and \eqref{C32}, we have, for
$i>\max\{\vartheta_t^{-1},\tilde\delta{}^{-1}\}$, \begin{align}&
p(t-\myfrac1i,x_k,y)\,2\tau(\myfrac1i,\max\{\myfrac14,\myfrac18\,d(x_0,y)\})\notag\\
&\qquad\le C_t\exp\big[-d(x_k,y)^{1+\eps_t}\,\big]2
C\exp\big[-\max\{\myfrac14,\myfrac18\,d(x_0,y)\}\,\big]\notag\\&\qquad\le
\operatorname{const}\,\exp\big[-\myfrac{1}{16}\,d(x_k,y)-\max\{\myfrac14,\myfrac18\,d(x_0,y)\}
\,\big]\notag\\
&\qquad\le\operatorname{const}\,\exp\big[-\myfrac1{16}\,d(x_0,x_k)\big]\,
\exp\big[\myfrac1{16}\,d(x_0,y)-
\max\{\myfrac14,\myfrac18\,d(x_0,y)\}\,\big]
.\label{763476}\end{align} Hence, by \eqref{klaw},
\eqref{sersre}--\eqref{763476}, \eqref{lpseaq}, and the dominated
convergence theorem, we get
\eqref{huizgurfrztrft}.\quad$\blacksquare$\vspace{2mm}

From \eqref{klaw}, \eqref{awuk}, and Lemmas~\ref{lemmaeins}
and~\ref{newlemma}, we get the central lemma of the proof:

\begin{lem}\label{lemmazwei} Let $x\in\widehat X^\N$ satisfy
\eqref{klaw}\rom. For $t>0$ and  $i\ge{\cal I}_t$, we set
\begin{gather*}{\Bbb A}_{t,i}:=\bigg[ \bigcup_{l={\cal
I}_t}^i\Big\{\,
\omega\in\Omega^\N:\sum_{k=1}^\infty\exp\big[-\myfrac1n\,
d(x_0,\omega_k(t-\myfrac1l))\,\big]<\infty\ \forall n\in\N,\\
d(\omega_k(t-\myfrac1l),\omega_k(s))\le\max\{1,\myfrac12\,
d(x_0,\omega_k(t-\myfrac1l))\}\ \forall s\in(t-\myfrac1l,t]\
\forall k\in\N\,\Big\}\bigg]\\\bigcap \Big\{\, \omega\in\Omega^\N:
\sum_{k=1}^\infty\exp\big[-\myfrac1n\,d(x_0,\omega_k(t))\,\big]<\infty\
\forall n\in\N,\\ d(\omega_k(t),\omega_k(s))\le\max\{1,\myfrac12\,
d(x_0,\omega_k(t))\ \forall s\in (t,t+\myfrac1i],\ \forall
k\in\N\,\Big\},\end{gather*} and for $t=0$ and $i\in\N$\rom, we
set
\begin{gather*}{\Bbb A}_{t,i}={\Bbb A}_{0,i}:=\big\{\,\omega\in\Omega^\N:
d(x_k,\omega_k(s))\le\max\{1,\myfrac12\,d(x_0,x_k)\}\ \forall s\in
(0,\myfrac1i] \ \forall k\in \N\,\big\}.\end{gather*} Then\rom,
$${\Bbb P}_{\bf x}\bigg(\bigcup_{i={\cal I}_t}^\infty{\Bbb
A}_{t,i}\bigg)= \lim_{i\to\infty}{\Bbb P}_{\bf x}({\Bbb
A}_{t,i})=1$$ for each $t\in\R_+$\rom.
\end{lem}

Let $D:=\R_+\cap{\Bbb Q}$ ($\Bbb Q$ denoting the set of rational
numbers), and $D=\{t_l\}_{l=1}^\infty$. We consider arbitrary
numbers $\eps_{lp}>0$, $l,p\in\N$, such that
$\sum_{l=1}^\infty\eps_{lp}<\infty$ for each $p\in\N$, and
\begin{equation}\label{yxvb}\lim_{p\to\infty}
\sum_{l=1}^\infty\eps_{lp}=0.\end{equation} For each $l,p\in\N$,
we choose an $i_{lp}\in\N$ such  that
\begin{equation}\label{poserdk}{\Bbb P}_{\bf x}\big( {\Bbb A}^{\mathrm
c}_{t_l,i_{lp}} \big)\le\eps_{lp}\,,\end{equation} which exists
due to Lemma~\ref{lemmazwei}. Then, we set
\begin{equation}\label{xrest}{\Bbb A}_p:=\bigcap_{l=1}^\infty{\Bbb
A}_{t_l,i_{lp}}.\end{equation} By \eqref{yxvb}, \eqref{poserdk},
and \eqref{xrest},
\begin{equation}\label{lkawvb}\lim_{p\to\infty}{\Bbb P}_{\bf x}({\Bbb
A}_p)=1.\end{equation}

For each $p\in\N$, let us consider the set
$$T_p:=\R_+\cap\bigg[\bigcup_{l=1}^\infty(t_l-i_{lp}^{-1},t_l+i_{lp}^{-1})\bigg].$$
Since the set
$\bigcup_{l=1}^\infty(t_l-i_{lp}^{-1},t_l+i_{lp}^{-1})$ is open in
$\R$, and since $T_p$ is dense in $\R_+$, we have
$T_p=\R_+\setminus T_p^{\mathrm c}$, where $ T^{\mathrm c}_p$ is
countable. We set $${\bf A}_p:={\Bbb A}_p\cap\bigg[\bigcap_{t\in
T_p^{\mathrm c} }\bigcup_{i={\cal I}_{t}}^{\infty}{\Bbb
A}_{t,i}\bigg].$$ By Lemma~\ref{lemmazwei}, we get
\begin{equation}\label{rebjh}{\Bbb P}_{\bf x}({\bf A}_p)={\Bbb
P}_{\bf x}({\Bbb A}_p).\end{equation}

Finally, we set \begin{equation}\label{23459834}{\bf
A}:=\bigcup_{p=1}^\infty{\bf A}_p.\end{equation}
 Therefore, by \eqref{lkawvb} and
\eqref{rebjh},  we get
\begin{equation}\label{swet}{\Bbb P}_{\bf x}({\bf A})=1.\end{equation}

\begin{lem}\label{lemmadrei} For any $\omega\in{\bf A}$ and $n\in\N$\rom, we
have $$\forall t\in\R_+:\quad
B_n(\omega(t)):=\sum_{k=1}^\infty\exp\big[-\myfrac1n\,
d(x_0,\omega_k(t))\,\big]<\infty ,$$ and moreover\rom, the mapping
 $$\R\ni t\mapsto B_n(\omega(t))\in\R$$ is continuous\rom.\end{lem}

\noindent{\it Proof}. We note that, if $B_{2n}(\omega(t))<\infty$
and if there exists an interval $(a,b)\subset\R_+$ such that
$t\in[a,b]$ and
$$d(\omega_k(s),\omega_k(t))\le\max\{1,\myfrac12\,d(x_0,\omega_k(t))\},\qquad
s\in(a,b),$$ then
 the series
$$\sum_{k=1}^\infty\exp\big[-\myfrac1n\,d(x_0,\omega_k(s))\,\big],\qquad
s\in(a,b), $$ are majorized by the convergent series
$$\sum_{k=1}^\infty
\exp\big[-\myfrac1{2n}\,d(x_0,\omega_k(t))\,\big].$$ Hence, the
statement follows from the construction of the set $\bf
A$.\quad$\blacksquare$\vspace{2mm}

Now, we define the action of the group $S_\infty$ on $\Omega^\N$
by
$$\sigma((\omega_k)_{k=1}^\infty):=(\omega_{\sigma(k)})_{k=1}^\infty,\qquad
\sigma\in S_\infty.$$ Evidently, the set $\bf A$ is invariant
under the action of $S_\infty$. We introduce the factor space
$\Omega^\N/S_\infty$ consisting of factor  classes
$[\omega_1,\omega_2,\dots]$. Analogously to \eqref{oizgdr} and
\eqref{awghzj}, we introduce then the mapping
$$\Omega^\N\ni\omega=(\omega_k)_{k=1}^\infty\mapsto {\bf
I}\omega=[\omega_1,\omega_2, \dots]\in\Omega^\N/S_\infty.$$

\begin{lem}\label{lemmavier}We have $${\bf I}{\bf A}\subset
C(\R_+;\dd\Gamma_\infty),$$ where  $C(\R_+; \dd\Gamma_\infty)$
denotes the set of continuous mappings from $\R_+$ into
$\dd\Gamma_\infty.$\end{lem}

\noindent{\it Proof}. Taking notice of Lemma~\ref{lemmadrei} and
of the definition of the metric space $\dd\Gamma_\infty$, it
remains only to show that, for each fixed
$(\omega_k)_{k=1}^\infty\in\bf A$, the mapping
\begin{equation}\label{awbh}\R_+\ni
t\mapsto\{\omega_k(t)\}_{k=1}^\infty\in\dd\Gamma_\infty\end{equation}
is vaguely continuous.

To this end, let us fix any $t\in\R_+$ and any ball $B(x_0,r)$
 of radius $r>0$. Then, there exist
$\eps>0$ and $K\in\N$ such that $$\sum_{k=K}^\infty\exp\big[
-d(x_0,\omega_k(s))\,\big]<e^{-r},\qquad
s\in\R_+\cap(t-\eps,t+\eps)$$ (see the proof of
Lemma~\ref{lemmadrei}). Hence, $$\omega_k(s)\not\in
B(x_0,r),\qquad k\ge K,\ s\in\R_+\cap(t-\eps,t+\eps),$$ which,
together with the continuity of each $\omega_k$ as a mapping from
$\R_+$ into $X$, implies the vague continuity of
\eqref{awbh}.\quad $\blacksquare$\vspace{2mm}

Thus, by Lemma~\ref{lemmavier}, we have that $${\bf I}\colon {\bf
A}\to \pmb\Omega:=C(\R_+;\dd\Gamma_\infty). $$ Denote  the trace
$\sigma$-algebra of ${\cal C}_\sigma(\Omega^\N)$ on $\bf A$ by
${\cal C}_\sigma({\bf A})$. Let ${\bf F}$ be the product
$\sigma$-algebra on $\pmb \Omega=C(\R_+;\dd\Gamma_\infty)$
generated by the $\sigma$-algebra ${\cal B}(\dd\Gamma_\infty)$:
$${\bf F}:=\sigma\{{\bf X}_t,\,t\in\R_+\},$$ where $${\bf
X}_t(\omega):=\omega(t).$$  Since the mapping $I\colon \widehat
X^\N\to\dd\Gamma_X$ defined by \eqref{awghzj} is measurable and
since ${\cal B}(\dd\Gamma_\infty)$ is the trace $\sigma$-algebra
of ${\cal B}(\dd\Gamma_X)$ on $\dd\Gamma_\infty$, we easily
conclude that the mapping $\bf I$ is ${\cal C}_\sigma({\bf
A})$-$\bf F $-measurable. Because of \eqref{swet}, we can consider
${\Bbb P}_{\bf x}$ as a probability measure on $({\bf A},{\cal
C}_\sigma({\bf A}))$ and let ${\bf P}_{\bf x}$ denote the image of
this measure under the mapping $\bf I$. Thus, ${\bf P}_{\bf x}$ is
a probability measure on $(\pmb\Omega,{\bf F})$.

For each $\sigma\in S_\infty$, the measures ${\bf P}_{\bf x}$ and
${\bf P}_{\sigma(\bf x)}$ evidently coincide, and so for each
$\gamma\in\dd\Gamma_\infty$, we can introduce the probability
measure ${\bf P}_\gamma:={\bf P}_{\bf x}$, where $\bf x$ is an
arbitrary element of the set $I^{-1}\gamma$.

Finally, we introduce  the sub-$\sigma$-algebras ${\bf
F}_t:=\sigma\{{\bf X}_s,\, s\le t\}$ and the translations
$({\pmb\theta}_t\omega)(s):=\omega(s+t)$, $t\in\R_+$. Thus, we get
\begin{equation}\label{qagf}({\pmb \Omega},{\bf F},({\bf
F}_t)_{t\in\R_+},({\pmb\theta}_t)_{t\in\R_+},({\bf
P}_\gamma)_{\gamma\in\dd\Gamma_\infty},({\bf
X}_t)_{t\in\R_+}).\end{equation} It follows directly from our
construction that \eqref{qagf} is a realization of the independent
infinite particle process,

 For a fixed
$\gamma\in\dd\Gamma_\infty$, the finite-dimensional distributions
of the process ${\bf X}_t$ under ${\bf P}_\gamma$ are given by
\begin{gather*}{\bf P}_\gamma({\bf X}_{t_1}\in A_1,{\bf X}_{t_2}\in A_2,
\dots,{\bf X}_{t_n}\in A_n)={\bf P}_{\bf x}({\bf X}_{t_1}\in
A_1,{\bf X}_{t_2}\in A_2,\dots,{\bf X}_{t_n}\in A_n)\\ ={\Bbb
P}_{\bf x}(\omega(t_1)\in I^{-1}A_1,\omega(t_2)\in
I^{-1}A_2,\dots,\omega(t_n)\in I^{-1}A_n)\\=\int_{I^{-1}A_1}{\Bbb
P} _{t_1}({\bf x},d{\bf x}_1)\int_{I^{-1}A_2}{\Bbb P}
_{t_2-t_1}({\bf x}_1,d{\bf x}_2)\dots\int_{I^{-1}A_n}{\Bbb
P}_{t_n-t_{n-1}}({\bf x }_{n-1},d{\bf x}_n)\\=\int_{A_1}{\bf P}
_{t_1}(\gamma,d \gamma_1)\int_{A_2}{\bf P
}_{t_2-t_1}(\gamma_1,d\gamma_2)\dots\int_{A_n}{\bf
P}_{t_n-t_{n-1}}(\gamma_{n-1},d\gamma_n),\\
0<t_1<t_2<\dots<t_n,\quad A_1,A_2,\dots,A_n\in{\cal
B}(\dd\Gamma_\infty)\end{gather*} where ${\Bbb P }_t({\bf
x}_i,d{\bf x}_j):={\Bbb P}_{t,{\bf x}_i}(d{\bf x}_j)$ and $\bf x$
is an arbitrary element of $I^{-1}\{\gamma\}$. Thus, the
finite-dimensional distributions of ${\bf X}_t$ are determined by
the Markov semigroup of kernels $({\bf P}_t)_{t\in\R_+}$. Hence,
it follows  that \eqref{qagf} is a time homogeneous Markov process
on $(\dd\Gamma_\infty,{\cal B}(\dd\Gamma_\infty))$ with transition
probability function $({\bf P }_t)_{t\in\R_+}$ (see e.g.\
\cite[Ch.~1, Sect.~3]{Bl}).

Finally, we note that any measure on the space $({\pmb
\Omega},{\bf F})$ is uniquely  determined by its
finite-dimensional distributions, and therefore the constructed
continuous Markov process is unique. \quad$\blacksquare$

\begin{rem}\rom{It is easy to see that the process $({\bf
X}_t)_{t\in\R_+}$constructed in the course of the proof of
Theorem~\ref{klhjuzip} is even Markov with respect to the
filtration $({\bf F}_{t+})_{t\in\R_+}$, where ${\bf
F}_{t+}:=\bigcap_{s>t}{\bf F}_t$. }\end{rem}

The following corollary states that, if the dimension $d$ of the
manifold $X$ is $\ge2$, then the process ${\bf X}_t$ starting at
$\gamma\in\Gamma_\infty$ lives with $\P_\gamma$-probability one in
$\Gamma_\infty$, i.e., the particles never collide (compare with
\cite{RS}).

\begin{cor}\label{diagonal}
 Let \rom{(C$1'$), (C2),} and \rom{(C3)} hold\rom, and let $d\ge2$\rom.
Then\rom, the independent infinite particle process can be
realized as the unique  continuous\rom, time homogeneous Markov
process $${\bf M}=({\pmb{ \Omega}},{\bf F},({\bf
F}_t)_{t\in\R_+},({\pmb \theta}_t)_{t\in\R_+}, ({\bf P
}_\gamma)_{\gamma\in\Gamma_\infty},({\bf X}_t)_{t\in\R_+})$$ on
the state-space $(\Gamma_\infty,{\cal B}(\Gamma_\infty))$ with
transition probability function $(\PP_t)_{t\in\R_+}$\rom.
\end{cor}

\noindent{\it Proof}. First, we claim that, if $d\ge2$, then
\begin{equation}\label{lkdfse} P_{x_1}\otimes P_{x_2}\big(\exists
t>0:\omega_1(t)=\omega_2(t)\big)=0,\qquad x_1,x_2\in X,
\end{equation} i.e., two independent Brownian motions on $X$ never
collide.

In the Euclidean case $X=\R^d$, this is a direct consequence of a
classical result from potential theory. Indeed, $\omega_1(\myfrac
t2)-\omega_2(\myfrac t2)$ is standard Brownian motion on $\R^d$
starting at $x_1-x_2$, and therefore \eqref{lkdfse} is equivalent
to the equality $$P_{x_1-x_2}\big(\exists t>0:
\omega(t)=0\big)=0,$$ which is true since points are polar for
Brownian motion on $\R^d$ if $d\ge2$ (see e.g.\
\cite[Proposition~2.5]{PortStone}).

In the general case, to prove \eqref{lkdfse} one can follow the
idea of  \cite{RS}.  First, we note that \eqref{lkdfse} is
equivalent to
 \begin{equation}\label{lhewfse} P_{x_1}\otimes
P_{x_2}\big(\exists
t>\myfrac1n:\omega_1(t)=\omega_2(t)\big)=0\qquad\forall n\in\N.
\end{equation}
Using the Markov property, we have
\begin{gather}P_{x_1}\otimes P_{x_2}\big(\exists
t>\myfrac1n:\omega_1(t)=\omega_2(t)\big)=\notag\\
=\int_{X^2}m(dy_1)\,m(dy_2)\,
p(\myfrac1n,x_1,y_1)p(\myfrac1n,x_2,y_2)\,
 P_{y_1}\otimes
P_{y_2}\big(\exists
t>0:\omega_1(t)=\omega_2(t)\big).\label{9786}\end{gather} By
virtue of \eqref{lhewfse} and \eqref{9786}, it suffices to verify
that the equality \eqref{lkdfse} holds only for
$m^{\otimes2}$-a.a.\ $(x_1,x_2)\in X^2$.

There exists a countable, locally finite  covering
$\{U^{(i)}\}_{i=1}^\infty$ of the manifold $X$ such that each
$U^{(i)}$ is an open set in $X$ diffeomorphic to the open cube
$(-1,1)^d$ in $\R^d$. Furthermore, two independent Brownian
motions
 on $X$ which start respectively at $x_1$
and $x_2$ form a Brownian motion on the manifold $X^2$ starting at
the point $(x_1,x_2)$. Hence, our problem can be reduced to the
following one: Show that
\begin{gather}
\label{hjghawse}P_{(x_1,x_2)}\big(\exists t>0:
\omega_1(t)=\omega_2(t)\in U^{(i)}_j\big)=0,\qquad \text{for
$m^{\otimes2}$-a.a.\ $(x_1,x_2)\in X^2$, $i,j\in\N$},\end{gather}
where $U^{(i)}_j$ is the subset of $U^{(i)}$ that is diffeomorphic
to the open cube $${\cal C }_j:=(-1+(1+j)^{-1},1-(1+j)^{-1})^d.$$

Let us consider the Dirichlet form that corresponds to Brownian
motion on $X^2$: \begin{align}{\cal E}(f,g):=&\int_{X^2}\big [
\langle\nabla^X_{x_1}f(x_1,x_2),\nabla^X
_{x_1}g(x_1,x_2)\rangle_{T_{x_1}(X)}\notag\\&\qquad+\langle\nabla^X_{x_2}f(x_1,x_2),\nabla^X
_{x_2}g(x_1,x_2)\rangle_{T_{x_2}(X)}
\big]\,m(dx_1)\,m(dx_2).\label{Dirin}\end{align} The bilinear form
$\cal E$ is defined first for $f,g\in{\cal
D}^{\otimes2}=C_0^\infty(X^2)$, and then it is closed.

Since $m^{\otimes 2}\big(\{(x_1,x_2)\in X^2: x_1=x_2\}\big)=0$, to
prove \eqref{hjghawse} it is enough to construct a sequence
$\{u_n\}_{n=1}^\infty\subset\operatorname{Dom}({\cal E})$ such
that $u_n$'s converge pointwisely to the indicator function of the
set $\{(x_1,x_2)\in X^2: x_1=x_2\in U_j^{(i)}\}$ and $\sup_n{\cal
E}(u_n,u_n)<\infty$ (see \cite{RS}).

By using the representation of the Dirichlet form $\cal E$ in
local coordinates on $U^{(i)}$, we get, for any function
$f\in\operatorname{Dom}({\cal E})$ having support in
$(U^{(i)})^2$: \begin{multline}{\cal
E}(f,f)=\frac12\int_{(-1,1)^{2d}}\sum_{k,l=1}^d \bigg[
g^{kl}(x_1)\,\frac{\partial f}{\partial x_1^k}(x_1,x_2)\,\frac
{\partial f}{\partial x_1^l}(x_1,x_2)\\ \text{} +
g^{kl}(x_2)\,\frac{\partial f}{\partial x_2^k}(x_1,x_2)\,\frac
{\partial f}{\partial x_2^l}(x_1,x_2)\bigg] \,\sqrt{{\bf
g}(x_1)}\,\sqrt{{\bf g}(x_2)}\,dx_1^1\dotsm dx_1^d\,dx_2^1\dotsm
dx_2^d,\label{qwewgh}\end{multline} where $\bf g$ denotes the
determinant of the matrix $(g_{kl})_{k,l=1}^d:=(\langle
\frac{\partial}{\partial x^k},\frac{\partial}{\partial
x^l}\rangle)_{k,l=1}^d$, and $(g^{kl})_{k,l=1}^d$ is its inverse.
We conclude from \eqref{qwewgh} that there exists a constant $C>0$
such that, for each function $f\in\operatorname{Dom}(\cal E)$
having support in a fixed $(U^{(i)}_j)^2$, $${\cal E}(f,f)\le C
{\cal E}_{\mathrm E}(f,f),$$ where $ {\cal E}_{\mathrm E}$ is the
(Euclidean) Dirichlet form on $\R^{2d}$: $${\cal E}_{\mathrm
E}(f,f)= \int_{\R^{2d}}\langle\nabla f(x_1,x_2),\nabla
f(x_1,x_2)\rangle\,dx_1^1\dotsm dx_1^d\,dx_2^1\dotsm dx_2^d,$$
$\nabla$ denoting the usual gradient on $\R^{2d}$. Hence, it
suffices to construct for any fixed $j\in\N$ a sequence
$\{u_n\}_{n=1}^\infty\subset\operatorname{Dom}({\cal E}_{\mathrm E
})$ such that each $u_n$ has support in ${\cal C}_{j+1}^2$, the
$u_n$'s converge pointwisely to the indicator of the set
$\{(x_1,x_2)\in\R^{2d}:x_1=x_2\in{\cal C}_j\}$, and $\sup_n{\cal
E}_{\mathrm E}(u_n,u_n)<\infty$. But the existence of such
sequence can be seen by a trivial modification of the proof of
Proposition~1 in \cite{RS}. Thus, \eqref{hjghawse} and hence also
the claim \eqref{lkdfse} are proven.

The rest of the proof follows from that of Theorem~\ref{klhjuzip}.
Instead of the set $\bf A$ given by \eqref{23459834}, one should
use its subset
\begin{equation}\label{893276}{\bf A}':={\bf A}\cap
 \bigg[\bigcap_{\{i,j\}\subset\N}\big\{\,\omega_i(t)\ne\omega_j(t)\
\forall t\in\R_+ \,\big\} \bigg].\end{equation} By \eqref{oiwegf},
\eqref{swet}, \eqref{lkdfse}, and \eqref{893276}, we get $${\Bbb
P}_{\bf x}({\bf A}')=1$$ for each ${\bf x}\in \tilde X^\N$
satisfying \eqref{klaw}, so that the measure ${\Bbb P}_{\bf x}$
can be considered as a probability measure on $({\bf A}',{\cal
B}({\bf A}'))$. Finally, noting that $${\bf I}{\bf A}'\subset
C(\R_+;\Gamma_\infty)$$ (compare with Lemma~\ref{lemmavier}), we
get the corollary by a corresponding modification of the last part
of the proof of Theorem~\ref{klhjuzip}.\quad$\blacksquare$

\begin{rem}\rom{(Path-wise construction of the independent particle process
 on $\Gamma_\infty$)  Let $d\ge2$ and let us consider the
probability space $(\pmb\Omega,{\cal C}_\sigma(\pmb\Omega),{\bf
P})$, where $\pmb\Omega:=\Omega^X$, $\Omega:=C(\R_+;X)$, ${\cal
C}_\sigma(\pmb \Omega)$ is the product $\sigma$-algebra on $\pmb
\Omega$ constructed from the $\sigma$-algebra $\cal F$ on $\Omega$
defined by \eqref{347845},  and  ${\bf P }:=\bigotimes_{x\in
X}P_x$. For any fixed $\gamma\in\Gamma_\infty$ and $t\in\R_+$, we
define $$\pmb \Omega\ni\omega=(\omega_x)_{x\in X}\mapsto {\bf
X}^\gamma_t(\omega):=\sum_{x\in\gamma}\eps_{X_t^x(\omega)},$$
where $X_t^x(\omega):=\omega_x(t)$. Thus, for any
$\gamma\in\Gamma_\infty$ we have constructed a process ${\bf
X}^\gamma:=({\bf X}_t^\gamma)_{t\in\R_+}$ which takes values in
the space of all measures on $X$. Let us fix any ${\bf
x}=(x_k)_{k=1}^\infty\in I^{-1}\{\gamma\}$. Then, ${\Bbb P}_{\bf x
}={\bf P}\circ {\bf I}_{\bf x}^{-1}$, where ${\Bbb P}_{\bf x }$ is
defined by \eqref{oiwegf}
 and $$\pmb \Omega\ni\omega=(\omega_x)_{x\in
X}\mapsto {\bf I}_{\bf
x}\omega:=(\omega_{x_k})_{k=1}^\infty\in\Omega^\N.$$ Hence, it
follows from the proof of Theorem~\ref{klhjuzip} (respectively
Corollary~\ref{diagonal}) that with $\bf P$ probability one the
independent particle process ${\bf X}^\gamma$ starts at $\gamma$,
never leaves $\Gamma_\infty$, i.e., ${\bf P}\big(\forall t>0: {\bf
X}^\gamma_t\in\Gamma_\infty\big)=1$, and has sample paths which
are continuous in the $d_\infty$ metric.}

\rom{In the case $d=1$, in order to give a corresponding path-wise
construction of the independent particle process on
$\dd\Gamma_\infty$, we proceed as follows. We consider the
probability space $({\pmb{\Omega}},{\cal
C}_\sigma({\pmb{\Omega}}),{\bf P})$, where $$
{\pmb{\Omega}}:=\underset{n(\cdot)\in\N^X}{\times}\underset{x\in
X}{\times}\Omega^{n(x)},$$ ${\cal C}_\sigma({\pmb\Omega})$ is the
corresponding product $\sigma$-algebra on ${\pmb\Omega}$, and
$${\bf P}:=\bigotimes_{n(\cdot)\in\N^X}\bigotimes_{x\in
X}P_x^{\otimes n(x) }.$$ For any $\gamma\in\dd\Gamma_X$, we define
$\hat\gamma\in\Gamma_X$ and a mapping $n_\gamma:\hat\gamma\to\N$
by $$\hat\gamma:=\operatorname{supp}\gamma,\qquad \hat\gamma\ni
x\mapsto n_\gamma(x):=\gamma(\{x\})\in\N.$$  We extend the mapping
$n_\gamma(\cdot)$ to the whole of $X$ by $n_\gamma(x):=1$ for all
$x\in X\setminus\hat\gamma$. Now, for any fixed
$\gamma\in\dd\Gamma_\infty$, we define a process ${\bf
X}^\gamma:=({\bf X}_t^\gamma)_{t\in\R_+}$ setting for each
$t\in\R_+$ $${\pmb\Omega}\ni\omega=(\omega^1_{n(\cdot),\,
x},\dots,\omega^{n(x)}_{n(\cdot),\,x })_{n(\cdot)\in\N^X,\, x\in
X}\mapsto  {\bf
X}^\gamma_t(\omega):=\sum_{x\in\hat\gamma}\sum_{i=1}^{n_\gamma(x)}\eps_{\omega^i_{n_\gamma
(\cdot),\,x}(t)}.$$ Analogously to the above, we conclude that
with ${\bf P}$ probability one the independent particle process $
{\bf X}^\gamma$ starts at $\gamma$, never leaves
$\dd\Gamma_\infty$, and has sample paths continuous in the
$d_\infty$ metric. }\end{rem}

\begin{rem}\rom{
The Markov process on the state space $\Gamma_\infty$ that was
constructed  in Corollary~\ref{diagonal} is a strong Markov
process. This can be shown by a modification of the proof of
  \cite[Theorem~5.10]{dynkin} using Corollary~\ref{awnmtezd}.
Furthermore, by proving a corresponding Feller property of ${\bf P
}_t$ on $\dd\Gamma_\infty$ with respect to the metric $d_1$, one
can show that the Markov  process on the state space
$\dd\Gamma_\infty$ that was constructed in Theorem~\ref{klhjuzip}
also possesses the strong Markov property. }\end{rem}

\section{Appendix: Proof of Lemma~\ref{alaNelson}}

Let us fix $r>0$, $\delta>0$, $x\in X$, $0\le t_1<\dots<t_n$,
$n\ge2$, with $t_n-t_1\le\delta$. Let
\begin{equation}\label{setfr}A:=\big\{\,\om\in\Omega: d(\om(t_1),\om(t_j))>r\ \text{for
some }j=2,\dots,n\,\big\},\end{equation} and let us show that
\begin{equation}
\label{wasesese} P_x(A)\le
2\tau(\delta,\myfrac12\,r).\end{equation}

Let \begin{align*}B:&=\big\{\, \om\in\Omega:
d(\om(t_1),\om(t_n))>\myfrac12\,r\,\big\},\\ C_j:&=\big\{\,
\om\in\Omega: d(\om(t_j),\om(t_n))>\myfrac12\,r\,\big\},\\
D_j:&=\big\{\, \om\in\Omega: d(\om(t_1),\om(t_j))>r,\ \text{and
}d(\om(t_1),\om(t_k))\le r\ \text{for
}k=1,\dots,j-1\,\big\}.\end{align*} Then, $$A\subset
B\cup\bigg[\bigcup_{j=1}^n(C_j\cap D_j)\bigg].$$  Therefore,
$$P_x(A)\le P_x(B)+\sum_{j=1}^n P_x(C_j\cap D_j).$$

Define $\tilde D_j\subset X^j$, $\tilde C_j\subset X^2$ by
\begin{align*}
\tilde D_j:&=\big\{\,(x_1,\dots,x_j)\in X^j: d(x_1,x_j)>r\text{
and }d(x_1,x_k)\le r\text{ for }k=1,\dots,j-1\,\big\},\\ \tilde
C_j:&=\big\{\, (x_1,x_2) \in X^2: d(x_1,x_2)>\myfrac12\,r\,\big\}.
\end{align*}
Then,
\begin{align*}P_x(C_j\cap D_j)&=\int_X \dots\int_X
p(t_1,x,dx_1)\dotsm p(t_j-t_{j-1},x_{j-1},dx_j)\\ &\quad\times
p(t_n-t_j,x_j,dx_n)1_{\tilde D_j}(x_1,\dots,x_j)1_{\tilde
C_j}(x_j,x_n)\\&\le \tau(\delta,\myfrac12\, r)\int_X\dots\int_X
p(t_1,x,dx_1)\dotsm p(t_j-t_{j-1},x_{j-1},dx_j)1_{\tilde
D_j}(x_1,\dots,x_j)\\&=\tau(\delta,\myfrac12\,
r)P_x(D_j).\end{align*} Since the sets $D_j$ are disjoint, we have
$$\sum_{j=1}^n P_x(C_j\cap D_j)\le \sum_{j=1}^n \tau
(\delta,\myfrac12\,r)P_x(D_j)\le \tau (\delta,\myfrac12\,r), $$
and since $P_x(B)\le \tau(\delta,\myfrac12\,r)$, we get
\eqref{wasesese}.

It follows from \eqref{setfr}, \eqref{wasesese} that
\begin{equation}\label{wadrdrhj}P_x(d(\om(t_j),\om(t_k))>2r\ \text{for some }j,k,\
1\le j,k\le n)\le 2\tau(\delta,\myfrac12\,r).\end{equation}
Indeed, if $d(\om(t_j),\om(t_k))>2r$, then
$d(\om(t_1),\om(t_j))>r$ or $d(\om(t_1),\om(t_k))>r$.

Since the estimate \eqref{wadrdrhj} is independent of $n$ and
$t_1,\dots,t_n$, we get
\begin{equation}\label{awfeshujr}P_x(d(\om(t_1),\om(t_2))>2r\ \text{for some
}t_1,t_2\in\R_+\cap{\Bbb Q},\ |t_1-t_2|\le \delta)\le 2
\tau(\delta,\myfrac12\,r).\end{equation}  Due to the continuity of
the trajectories $\om\in\Omega$, \eqref{awfeshujr} implies the
statement of the lemma.\quad $\blacksquare$

\begin{center}

{\bf Acknowledgments}
\end{center}

\noindent We thank Tobias Kuna for useful discussions. The authors
were partially supported by the SFB 256, Bonn University, and SFB
343, Bielefeld University. The support of the DFG through
Projects~436~113/39 and 436~113/43, and of the BMBF through
Project~UKR-004-99 is  gratefully
 acknowledged.


\end{document}